\documentclass[journal,11pt,twocolumn]{IEEEtran}
\usepackage{amsmath}
\usepackage{amssymb}
\usepackage{amsthm}
\usepackage{amsfonts}
\usepackage{amsxtra}
\usepackage{mathabx}
\usepackage{algorithm}
\usepackage{algorithmic}
\usepackage[pdfencoding=auto]{hyperref}
\usepackage[stable]{footmisc}
\usepackage{cleveref}
\usepackage{bookmark}
\usepackage{tikz-cd}
\usepackage{subfigure}
\usepackage{graphicx}
\usepackage{url}
\usepackage{booktabs}
\usepackage{nicefrac}
\usepackage{microtype}
\usepackage{mathtools}
\usepackage{bm}
\usepackage{bbm}
\usepackage{epstopdf}
\usepackage{multirow}
\usepackage[numbers,sort&compress]{natbib}
\usepackage[normalem]{ulem}

\newcommand{\R}{\mathbb{R}}

\renewcommand{\P}[1]{\operatorname{\mathbb{P}}\left(#1\right)}
\newcommand{\E}{\operatorname{\mathbb{E}}}

\newcommand{\vct}[1]{\mathbf{#1}}
\newcommand{\norm}[1]{\left \|#1 \right \|}

\newcommand{\mtx}[1]{\mathbf{#1}}

\newcommand{\trace}{\operatorname{trace}}

\newcommand{\calH}{\mathcal{H}}

\newcommand{\va}{\vct{a}}
\newcommand{\vb}{\vct{b}}

\newcommand{\ve}{\vct{e}}
\newcommand{\vf}{\vct{f}}

\newcommand{\vx}{\vct{x}}
\newcommand{\vy}{\vct{y}}

\newcommand{\vzero}{\vct{0}}
\newcommand{\vxi}{\vct{\xi}}

\newcommand{\mA}{\mtx{A}}
\newcommand{\mB}{\mtx{B}}

\newcommand{\mE}{\mtx{E}}

\newcommand{\mH}{\mtx{H}}
\newcommand{\mI}{\mtx{I}}

\newcommand{\mQ}{\mtx{Q}}

\newcommand{\mS}{\mtx{S}}
\newcommand{\mT}{\mtx{T}}
\newcommand{\mU}{\mtx{U}}
\newcommand{\mV}{\mtx{V}}
\newcommand{\mW}{\mtx{W}}
\newcommand{\mX}{\mtx{X}}

\newcommand{\mZ}{\mtx{Z}}
\newcommand{\mPhi}{\mtx{\Phi}}
\newcommand{\mPsi}{\mtx{\Psi}}

 \newtheorem{lemma}{Lemma}
 \newtheorem{theorem}{Theorem}
 
 \newtheorem*{remark}{Remark}

\newcommand{\xspace}{\mkern-1.5mu}
\newcommand\tnorm[1]{\left\vert\xspace\left\vert\xspace\left\vert\mskip2mu
#1\mskip2mu \right\vert\xspace\right\vert\xspace\right\vert}

\begin{document}
\title{Sketching low-rank matrices with a shared column space by convex programming}
\author{Rakshith S Srinivasa, Seonho Kim, and Kiryung Lee,~\IEEEmembership{Senior Member,~IEEE}\thanks{This work was supported in part by NSF CAREER Award CCF 19-43201. A preliminary version was presented at the 2021 Asilomar Conference on Signals, Systems, and Computers.}
\thanks{Rakshith S. Srinivasa is with Samsung Research America, Mountain View, CA 94043. 
Seonho Kim and Kiryung Lee are with the Department of Electrical and Computer Engineering at the Ohio State University, Columbus, OH 43220 (e-mail: 
kiryung@ece.osu.edu).}
}


\maketitle

\begin{abstract}
In many practical applications including remote sensing, multi-task learning, and multi-spectrum imaging, data are described as a set of matrices sharing a common column space.
We consider the joint estimation of such matrices from their noisy linear measurements. 
We study a convex estimator regularized by a pair of matrix norms. 
The measurement model corresponds to block-wise sensing and the reconstruction is possible only when the total energy is well distributed over blocks. 
The first norm, which is the maximum-block-Frobenius norm, favors such a solution. 
This condition is analogous to the notion of low-spikiness in matrix completion or column-wise sensing. 
The second norm, which is a tensor norm on a pair of suitable Banach spaces, induces low-rankness in the solution together with the first norm. 
We demonstrate that the joint estimation provides a significant gain over the individual recovery of each matrix when the number of matrices sharing a column space and the ambient dimension of the shared column space are large relative to the number of columns in each matrix.
The convex estimator is cast as a semidefinite program and an efficient ADMM algorithm is derived. 
The empirical behavior of the convex estimator is illustrated using Monte Carlo simulations and recovery performance is compared to existing methods in the literature. 
\end{abstract}


\section{Introduction}

We consider the problem of joint reconstruction of rank-$r$ matrices $\mX_1,\dots,\mX_K \in \R^{M \times N}$, which share a common left factor $\mU \in \R^{M \times r}$ from the linear measurements expressed as
\begin{align}
    y_{l,k} & = \langle \mB_{l,k}, \mX_k \rangle + w_{l,k}, \label{eq:sensing_model2}
\end{align} 
with $l \in [L] : = \{1,\dots,L\}$ and $k \in [K] := \{1,\dots,K\}$ and where $\mB_{l,k}$'s are a set of known ``sensing matrices'' and $w_{l,k}$'s represent additive noise in the measurements. 
Due to the assumption that $\mX_1,\dots,\mX_K$ share a common left factor, there exist $\mV_1, \dots, \mV_K \in \R^{N \times r}$ such that $\mX_k = \mU \mV_k^*$ for all $k \in [K]$. 
Let $\mX = [\mX_1 ~ \mX_2 ~ \cdots ~ \mX_K] \in \R^{M \times NK}$. Then each $y_{l,k}$ corresponds to a linear measurement of $\mX$ given by 
\begin{align}
    y_{l,k} & = \langle \mA_{l,k}, \mX \rangle + w_{l,k}, \label{eq:sensing_model}
\end{align} 
where $\mA_{l,k} = \ve_{k}^* \otimes \mB_{l,k}$ for all $l \in [L]$ and $k \in [K]$. 
In other words, the joint reconstruction of $\mX_1,\dots,\mX_K$ is considered as ``block-wise'' sensing of the concatenated rank-$r$ matrix $\mX$. 

The inverse problem for the model in \eqref{eq:sensing_model} has been investigated as a shared low-rank matrix regression in \cite{gigi2020spectral}. The authors provided a solution to the problem by using non-convex optimization to regress the shared subspace and the individual right factors separately. In particular, the authors propose to use spectral initialization followed by covariance estimation to solve for the shared subspace. With this estimate, they further use ridge regression to estimate the right factors. 

This problem arises naturally in numerous practical applications including remote sensing using satellite data \cite{gigi2020spectral}, multi-task learning \cite{lozano2012multitask, zhang2021survey}, and multi-channel data acquisition \cite{aziz2016lowrank}. 
Other applications include data compression in scientific simulations and multi-spectrum imaging. 
For example, a similar sketching problem of linear dimensionality reduction of streaming data has been considered \cite{tropp2019streaming}, in which data generated during simulations of fluid dynamics are shown to have a low-rank structure. 
The sensing model in \eqref{eq:sensing_model} applies to this application in the following sense: blocks of data generated over time (in this case, the data corresponds to the state of fluid motion over time) can be sketched independently, since they share a low-rank structure. 
Yet in another example, in the \textit{Square Kilometer Array} (SKA) \cite{spencer2013square}, astronomical data are collected using antenna elements spreading across different continents. 
Astronomical data collected via multi-channel acquisition show a low-rank structure via a shared factor \cite{aziz2016lowrank}. A similar multi-channel sensing application was also considered in \cite{srinivasa2019imaging}, where the linear model for data acquisition at different frequencies follow a nested subspace structure. 
Hence, the model in \eqref{eq:sensing_model} is directly applicable to the problem of dimensionality reduction before data aggregation.

With a rich context of relevant applications, this paper investigates a fundamental question on the inverse problem in the two equivalent models in \eqref{eq:sensing_model2} and \eqref{eq:sensing_model}. 
The main question we address is whether one can obtain a significant gain from the joint estimation by taking advantage of the redundancy across the matrices.
We focus on statistical analysis in a scenario, where the measurement matrices are independent copies of a random matrix with i.i.d. Gaussian entries of zero mean and unit variance, and the additive noise terms are independent of the signal and i.i.d. Gaussian. 

Note that the block-wise measurement model in \eqref{eq:sensing_model} takes linear measurements from a selected block at a time. 
In an extreme case, where all blocks of $\mX$ are zero matrices except a single block of unknown index, the measurements from zero blocks do not carry any useful information about $\mX$. 
Since the index of the nonzero block is unknown, $\tilde{O}(r(M+N))$\footnote{The tilde-big-O notation is defined as follows: If $a = \tilde{O}(b)$, then $a$ is less than $b$ times a logarithmic factor of considered parameters.} measurements per block are needed for ``stable'' recovery . 
In other words, there is no gain from joint estimation. 

One expects that a gain is achieved when the total energy of $\mX$ is well distributed across all blocks $\mX_1,\dots,\mX_K$. 
To favor a solution with this property, we consider an estimator regularized by the ``maximum'' correlation of $\mX$ with all possible measurement matrices $\mA_{k,l}$'s. 
Let $\gamma_{l,k} := \langle \mA_{l,k}, \mX \rangle = \langle \mB_{l,k}, \mX_k \rangle$. 
Then $\gamma_{l,k}$'s are i.i.d. Gaussian. 
Recall that the maximum of i.i.d. Gaussian random variables is upper-bounded with high probability by the standard deviation within a logarithmic factor of the number of random variables.
Following this observation, we consider the maximum correlation represented by
\begin{equation}
\label{eq:max_cor}
\max_{k \in [K], l \in [L]} \sqrt{\mathbb{E} \gamma_{l,k}^2}
= \max_{k \in [K]} \norm{\mX_k}_\mathrm{F}.
\end{equation}
The right-hand side of \eqref{eq:max_cor} is called the \textit{maximum-block-Frobenius} norm and will be denoted by
\begin{equation}
\label{eq:infty_Freb_norm}
\begin{aligned}
\norm{[\mX_1 ~ \mX_2 ~ \cdots ~ \mX_K]}_{\infty,\mathrm{F}}
& = \norm{\sum_{k=1}^K \ve_k^* \otimes \mX_k}_{\infty,\mathrm{F}} \\
& = \max_{k \in [K]} \norm{\mX_k}_\mathrm{F},    
\end{aligned}
\end{equation}
where $\ve_k \in \R^K$ denotes the $k$th column of the $K$-by-$K$ identity matrix $\mI_K$ for $k \in [K]$. 

To account for the low-rankness of $\mX$, we introduce another regularizer by a matrix norm given by 
\begin{equation}
\label{eq:dollar_norm}
     \norm{\mX}_{\$} = \inf_{\mU,\mV: \mU\mV^* = \mX} \norm{\mU}_\mathrm{F} \|\mV^*\|_{\infty,\mathrm{F}},
\end{equation}
where the common number of columns in $\mU$ and $\mV$ can be arbitrary while their product $\mU \mV^*$ coincides $\mX$. 
In general, matrix norms are not necessarily easy to compute. 
However, $\norm{\mX}_\$$ can be computed via a standard semidefinite program.
Note that the optimization in \eqref{eq:dollar_norm} is equivalent to 
\begin{equation*}
 \norm{\mX}_{\$} = \inf_{\mU,\mV: \mU\mV^* = \mX} \max \left(\|\mU\|_\mathrm{F}^2, \|\mV^*\|_{\infty,\mathrm{F}}^2\right), 
\end{equation*}
where the infimum is achieved if $\|\mU\|_\mathrm{F}$ coincides with $\|\mV^*\|_{\infty,\mathrm{F}}$. 
Furthermore, it has been shown (e.g. \cite{srebro2004maximum}) that there exist $\mU$ and $\mV$ such that $\mX = \mU \mV^*$, $\mW_1 = \mU \mU^*$, and $\mW_2 = \mV \mV^*$ if and only if 
\[
\begin{bmatrix} \mW_1 & \mX \\ \mX^* & \mW_2 \end{bmatrix} \succeq \vzero.
\]
Then we have $\norm{\mU}_\mathrm{F}^2 = \mathrm{trace}(\mW_1)$ and $\norm{\mV^*}_{\infty,\mathrm{F}}^2 = \max_{k \in [K]} \mathrm{trace} ((\ve_k^* \otimes I_N) \mW_2 (\ve_k \otimes I_N) )$. 
Therefore, $\norm{\mX}_\$$ can be computed via the following program:
\begin{equation}
\label{eq:compute_dollarnorm}
\def\arraystretch{1.25}
\begin{array}{lcl}
\norm{\mX}_{\$} = & \displaystyle \mathop{\mathrm{min}}_{\beta, \mW_1, \mW_2} & \beta \\
& \mathrm{s.t.} & \displaystyle \mathrm{trace}(\mW_1) \leq \beta \\
& & \displaystyle 
\mathrm{trace}\left( (\ve_k^* \otimes \mI_N) \mW_2 (\ve_k \otimes \mI_N) \right) \leq \beta, \\ 
& & \hspace{4.35cm} \forall k \in [K] \\
& & \displaystyle \begin{bmatrix}
\mW_1 & \mX \\
\mX^* & \mW_2
\end{bmatrix} \succeq \mathbf{0}.    
\end{array}
\end{equation}

The following lemma, proved in Appendix~\ref{sec:proof:lem:interlacing}, demonstrates how the above two norms characterize low-rankness through interlacing inequalities.

\begin{lemma}
\label{lem:interlacing}
Suppose that $\mX \in \R^{M \times NK}$ satisfies $\mathrm{rank}(\mX) \leq r$. Then we have
\begin{equation}
\label{eq:interlacing}
\norm{\mX}_{\infty,\mathrm{F}} 
\leq 
\norm{\mX}_{\$} 
\leq 
\sqrt{r} \norm{\mX}_{\infty,\mathrm{F}}.
\end{equation}
\end{lemma} 

We consider an estimator given as the solution to the following optimization program that minimizes the quadratic loss constrained to the two inequality constraints given by the above norm regularizers: 
\begin{equation}
\label{eq:estimator}
\mathop{\mathrm{minimize}}_{\mX \in \kappa(\alpha, \beta)} ~ \sum_{l=1}^L \sum_{k=1}^K \left( y_{l,k} - \langle \ve_{k}^* \otimes  \mB_{l,k}, \mX \rangle  \right)^2,
\end{equation}
where the constraint set is given by
\begin{equation}
\label{eq:def_kappa_ab}
    \kappa(\alpha, \beta) = \{\mX \in \R^{M \times NK}: \norm{\mX}_{\infty,\mathrm{F}} \leq \alpha,~ \norm{\mX}_{\$} \leq \beta \}.
\end{equation} 
Due to the characterization of the $\$$-norm in \eqref{eq:compute_dollarnorm}, the convex estimator in \eqref{eq:estimator} is obtained as a solution to 
\begin{equation}
\label{eq:estimator_sdp}
\def\arraystretch{1.25}
\begin{array}{lcl}
& \displaystyle \mathop{\mathrm{minimize}}_{\mX, \mW_1, \mW_2} ~ & \displaystyle \sum_{l=1}^L \sum_{k=1}^K \left( y_{l,k} - \langle \mB_{l,k}, \mX (\ve_k \otimes \mI_N) \rangle  \right)^2 \\
& \mathrm{subject~to} &\displaystyle \mathrm{trace}(\mW_{1}) \leq \beta \\
& & \displaystyle 
\mathrm{trace}\left((\ve_k^* \otimes \mI_N) \mW_2 (\ve_k \otimes \mI_N)\right) \leq \beta, \\
& & \hspace{4.55cm} k \in [K], \\
& & \left\|\mX\right\|_{\infty,\mathrm{F}}\leq\alpha\\
& & \displaystyle \begin{bmatrix}
\mW_{1} & \mX \\
\mX^* & \mW_{2}
\end{bmatrix} \succeq \mathbf{0}. 
\end{array}
\end{equation}

Our main results characterize the estimation problem with respect to the model $\kappa(\alpha,\beta)$ by an achievable error bound and a minimax lower bound. 
We first present an upper bound on the estimation error by convex program in \eqref{eq:estimator} in the following theorem. 

\begin{theorem}
\label{thm:main}
Let $(y_{l,k})$ be measurements of blocks of $\mX \in \R^{M \times KN}$ as described in \eqref{eq:sensing_model}. 
Suppose that $\mB_{l,k}$'s are independent copies of a random matrix whose entries are drawn i.i.d. from $\mathcal{N}(0,1)$. 
Furthermore, suppose that the noise entries $\gamma_{l,n}$'s are drawn from $\mathcal{N}(0,\sigma^2)$ and independent from everything else.
Then there exists a numerical constant $C$ such that if 
\begin{equation}
\label{eq:sample_rate}
L \geq C (\beta/\alpha)^2 
N \left(N+\frac{M}{K}\right) (\ln K)^3, 
\end{equation}
then it holds with probability $1-\zeta$ that the estimate $\widehat{\mX}$ of $\mX$ by \eqref{eq:estimator} satisfies
\begin{equation}
\label{eq:err_bnd}
\begin{aligned}
\|\widehat{\mX} - \mX\|_\mathrm{F}^2 
& \lesssim 
K \alpha^2 \left( 1 \vee \frac{\sigma}{\alpha} \right) \\
& \cdot \sqrt{\frac{
(\beta/\alpha)^2 N(M+NK) (\ln K)^3 + \ln(1/\zeta)}{LK}}    
\end{aligned}
\end{equation}
for all $\mX \in \kappa(\alpha,\beta)$.\footnote{We use a shorthand notation for the minimum and maximum of two numbers given by $\min(a,b) = a \wedge b$ and $\max(a,b) = a \vee b$.}
\end{theorem}

To interpret the result of Theorem~\ref{thm:main} in the context of joint estimation, we introduce the \textit{spikiness} parameter $\mu$ defined by
\[
\mu := \frac{\sqrt{K} \norm{\mX}_{\infty,\mathrm{F}}}{\norm{\mX}_\mathrm{F}}.
\]
The parameter $\mu$ of $\mX$ represents how the total energy of $\mX$ spreads over the blocks. 
A larger $\mu$ implies that there exist few blocks consuming most of the total energy. 
We also define the signal-to-noise-ratio (SNR) by 
\[
\mathrm{SNR} = \frac{\sum_{l=1}^L \sum_{k=1}^K \mathbb{E}[\langle \mA_{l,k},  \mX \rangle^2]}{\sum_{l=1}^L \sum_{k=1}^K \mathbb{E}[w_{l,k}^2]}
= \frac{\norm{\mX}_\mathrm{F}^2}{K \sigma^2}.
\]
Then the error bound in \eqref{eq:err_bnd} is rewritten as
\begin{equation}
\label{eq:err_bnd2}
\begin{aligned}
& \frac{\|\widehat{\mX} - \mX\|_\mathrm{F}^2}{\norm{\mX}_\mathrm{F}^2} \\
& \lesssim 
\mu^2 \left( 1 \vee \frac{\mu^{-1}}{\mathrm{SNR}^{1/2}} \right) \sqrt{\frac{(\beta/\alpha)^2 N(M+NK) (\ln K)^3}{LK}}.
\end{aligned}
\end{equation}
Furthermore, in a low-SNR regime, where $\mathrm{SNR} = O(\mu^{-2} \ln(LK))$, the error bound in \eqref{eq:err_bnd2} reduces to 
\begin{equation*}
\frac{\|\widehat{\mX} - \mX\|_\mathrm{F}^2}{\norm{\mX}_\mathrm{F}^2}
\lesssim 
\frac{\mu}{\mathrm{SNR}^{1/2}} \sqrt{\frac{
(\beta/\alpha)^2 N(M+NK) (\ln K)^3}{LK}}.
\end{equation*}
This implies that for a fixed SNR, the error decays as $\tilde{O}\left(\sqrt{\frac{\mu^2 (\beta/\alpha)^2 N(M/K + N)}{L}}\right)$. 
Note that the spikiness parameter $\mu$ in the error bound remains the same regardless of the distribution of columns norms within each block.

The upper bound on the estimation error becomes tightest when $\alpha = \norm{\mX}_{\infty,\mathrm{F}}$ and $\beta = \norm{\mX}_\$$. 
In practice, one needs to estimate those parameters so that $\alpha$ and $\beta$ are no less than the corresponding norms of $\mX_0$. 
To illustrate the optimal performance, suppose that $\alpha = \norm{\mX}_{\infty,\mathrm{F}}$, $\beta = \norm{\mX}_\$$, and $\mathrm{rank}(\mX) \leq r$. Then, by Lemma~\ref{lem:interlacing}, we have 
\[
\beta = \norm{\mX}_\$ \leq \sqrt{r} \norm{\mX}_{\infty,\mathrm{F}} \leq \sqrt{r} \alpha,
\]
which implies $(\beta/\alpha)^2 \leq r$. 
In the current scenario, the individual recovery of each block can succeed from $\tilde{O}(r(M+N))$ samples per block, but the joint recovery succeeds with $\tilde{O}(rN(\frac{M}{K}+N))$ samples per block. 
Therefore, if $M > N^2$ and $K > N$, then the joint recovery is feasible from fewer observations than the individual recovery. The advantage of our method is more pronounced for larger $M$ and $K$ (relative to $N$). 
For example, in the context of regression on hyperspectral remote sensing data, $M$ and $K$ respectively counts spectral bands and temporal samples while $N$ measures the size of a neighborhood of the target location in pixels, from which the prediction is made. 
Typical hyperspectral instruments have more than $200$ spectral bands \cite{bioucas2013hyperspectral}. 
Furthermore, it is feasible to learn the regressor from a large number of temporal samples. 
In this illustration, the parameters $M$ and $K$ are large relative to $N$. Hence, as discussed above, the joint recovery shows a significant gain over the individual recovery.

Next, we compare the upper bound in Theorem~\ref{thm:main} to a matching minimax lower bound.

\begin{theorem}
\label{thm:minimax}
Suppose that $(\beta/\alpha)^2(M \vee NK) \geq 48$. Then the minimax $\norm{\cdot}_\mathrm{F}$-risk is lower-bounded as   
\begin{align*}
& {\inf_{\widehat{\mX}}} {\sup_{\mX\in \kappa(\alpha, \beta)}} 
\frac{1}{K} \E \|\widehat{\mX} - \mX\|_\mathrm{F}^2 \\
& \geq \frac{\alpha^2}{16} \left( 1 \wedge \frac{\sigma}{8\sqrt{2}\alpha} \sqrt{\frac{(\beta/\alpha)^2(M \vee NK)}{L K}} \right).
\end{align*}    
\end{theorem}

\noindent Compared to the minimax bound by Theorem~\ref{thm:minimax}, the error bound for the estimator of \eqref{eq:estimator} in Theorem ~\ref{thm:main} is sub-optimal in general. However, the bound is near-optimal when the noise factor dominates and $N=O(1)$.
The minimax error bound decays with a rate proportional to $1/\sqrt{L}$, which is slower than the optimal rate $\sim 1/L$. We suspect that this is due to the relaxation of the set of low-rank matrices to the convex set $\kappa(\alpha,\beta)$.
On the other hand, with the relaxed matrix model, it applies to matrices with modeling error, for example, to approximately low-rank matrices.

\noindent\textbf{Related prior results:} Recovery of low-rank matrices under a structured measurement model has been of interest for many years with various applications in signal processing and statistics \cite{davenport2016overview}. 
Our approach is aligned with how the matrix completion problem was tackled with nuclear norm \cite{negahban2011estimation} and max norm \cite{cai2016matrix} without imposing the incoherence via singular value decomposition. 
A highly related model is column-wise sketching, which is a special case of \eqref{eq:sensing_model} with $N = 1$. 
Recent work provided sample complexity estimates using convex estimators \cite{lee2021approximately}. 
When $N=1$, the equivalence between the $2$-summing norm and the projective norm has been shown when a factor in the tensor product is equipped with the $\ell_\infty$ norm \cite[Lemma~4.4]{lee2021approximately}. However, the constraint set in \eqref{eq:def_kappa_ab} is determined by tensor norms on the product of two Banach spaces, neither of which uses the $\ell_\infty$ norm. Hence, 
even though Theorem~\ref{thm:main} produces the analogous result for $N=1$ \cite[Theorem~1.2]{lee2021approximately}, the extension in the other direction is not trivial. 
Therefore, the scenario with $N >1$ considered in this paper is significantly different from the case when $N = 1$. 
Importantly, as discussed earlier, there are applications modelled only by $N>1$. 

To the best of our knowledge, there is only one paper which studied the exact inverse problem in \eqref{eq:sensing_model2}. 
It has been shown that the spectral method provides an $\varepsilon$-accurate estimate of the column space of $\mU$, where the error is measured by the sine of the largest principal angle, from $O(\frac{\epsilon^{-2} r^4 M}{K} + N)$ noise-free samples per block with high probability \cite{gigi2020spectral}. 
In this paper, we improve upon their work in the following aspects: 
First, they only considered the recovery of only the column space of $\mU$ instead of $\mU \mV^*$, whereas the convex estimator in \eqref{eq:estimator} recovers the entire matrix. 
Second, our analysis continues to hold in the presence of measurement noise and model error, unlike the analysis in \cite{gigi2020spectral} which expects noise-free measurements.
Third, the unknown matrix in their analysis is arbitrarily fixed. Therefore, the error probability $O(\frac{1}{M})$ increases proportionally to the number of instances as one repeatedly applies the error bound to multiple instances. 
On the contrary, the error bound by Theorem~\ref{thm:main} provides a strong uniform guarantee that applies to all instances within the given model with high probability. 
It was proposed to further refine the estimate from the spectral method via gradient descent \cite{gigi2020spectral}. They demonstrated that the estimate by gradient descent from the spectral method outperforms that by random initialized gradient descent. 
In Section~\ref{sec:numerical}, we observed that gradient descent outperforms the convex estimator in \eqref{eq:estimator}. 
However, any error bound for the gradient descent estimator has not been established yet.

There has been a line of research on estimating low-rank matrices from structured measurements by iterative algorithms \cite{vaswani2017low,nayer2019phaseless,nayer2021sample,nayer2021fast}.
It has been shown that the ``sample-split'' version of alternating minimization and gradient descent from spectral initialization provides an $\epsilon$-accurate estimate from $\tilde{O}(r^2(M+K) \ln(1/\epsilon))$ noise-free phaseless measurements when the unknown matrix of size $M \times K$ is exactly rank-$r$. However, in practice, the sample-split algorithms perform significantly worse than the original counterpart. 
On the other hand, it has been shown that the vanilla gradient descent without sample splitting succeeds at a near optimal rate for phase retrieval, matrix completion, and blind deconvolution \cite{ma2020implicit}. 
However, it remains an open question whether the elegant analysis based on leave-one-out auxiliary sequences for gradient descent extends to the linear column-wise sensing. 
There also exists a convex optimization approach to low-rank recovery from phaseless measurements \cite{lee2021phase}. The considered linear models are different from the column-wise sensing but they have shown a near-optimal sample-complexity result without requiring sample splitting. 

The rest of this paper is organized as follows. Section~\ref{sec:definitions} introduces notation and definitions. 
Section~\ref{sec:entropy} derives the entropy estimate with respect to the $\$$-norm through its relation to the projective norm. The proof of Theorem~\ref{thm:main} is provided in Section~\ref{sec:proof:thm:main}, followed by discussions on numerical results in Section~\ref{sec:numerical}. 
We conclude with remarks and future directions in Section~\ref{sec:conclusion}.

\section{Notation}
\label{sec:definitions}

In this section, we introduce notation and definitions used throughout. 
Symbols for column vectors (resp. matrices) are denoted by boldface lower-case (upper-case) letters. 
For linear operator $\mT$ between vector spaces, the adjoint will be denoted by $\mT^*$. In a special case when $\mT$ is a matrix, then $\mT^*$ denotes the transpose. 
For vector space $X$, its algebraic dual is denoted by $X^*$. 
For Banach space $X$, the norm dual is denoted by $X^*$. 
The Kronecker product of two matrices $\mA$ and $\mB$ will be written as $\mA \otimes \mB$. 
The same symbol $\otimes$ is also used for general tensor product. 
We use various norms on column vectors and matrices throughout the paper. 
For column vector $\vx$, the $\ell_p$-norm is denoted by $\norm{\vx}_p$ for $p \geq 1$. 
Then the Banach space of column vectors of length $N$ with the $\ell_p$-norm is denoted by $\ell_p^N$. 
For matrix $\mA$, the Frobenius and spectral norms are denoted respectively by $\norm{\mA}_\mathrm{F}$ and $\norm{\mA}$. 
The corresponding unit norm balls are denoted by $B_\mathrm{F}$ and $B_\mathrm{S}$. 
More generally, the unit ball in a Banach space $X$ will be denoted by $B_X$. 
Furthermore, the operator norm of linear operator $\mT$ is written as $\norm{\mT}_\mathrm{op}$. 
For matrix $\mA \in \R^{M \times N}$, the column vector of length $MN$ obtained by stacking the columns of $\mA$ is denoted by $\mathrm{vec}(\mA)$.  
The maximum and minimum of two real numbers $a$ and $b$ will be respectively denoted by $a \vee b$ and $a \wedge b$. 

The convex estimator in \eqref{eq:estimator} induces a low-rank solution via the constraint set defined as in \eqref{eq:def_kappa_ab} by the max-block-Frobenius norm in \eqref{eq:infty_Freb_norm} and the $\$$-norm in \eqref{eq:dollar_norm}. 
The error analysis of the estimator is based on various properties of the $\$$-norm, which are characterized by tensor norms. 
A brief review on related mathematical background is provided in a companion paper \cite[Section~2]{lee2021approximately}. 
Further details can be found in monographs on tensor product \cite{defant1992tensor,diestel2008metric}.
Here we recall the minimal set of definitions which are necessary to state and derive the main results.

For vector spaces $X$ and $Y$, let $X^*$ and $Y^*$ denote the corresponding algebraic dual spaces, i.e. the collection of all linear functionals. 
The algebraic tensor product, denoted by $X \otimes Y$, is the set of all blinear functions on $X^* \times Y^*$. 
The algebraic tensor product is embedded into the set of all linear maps from $X^*$ to $Y$, denoted by $L(X^*,Y)$. 
In particular, if all vector spaces are finite dimensional, then $X \otimes Y$ is identified to $L(X^*,Y)$. 

Let $X$ and $Y$ be finite-dimensional Banach spaces. 
A norm on $X \otimes Y$ is a tensor norm if it satisfies 
\[
\norm{\vx \otimes \vy} \leq \norm{\vx}_X \norm{\vy}_Y, \quad \forall \vx \in X, ~ \vy \in Y
\]
and its dual norm satisfies 
\[
\norm{\vx^* \otimes \vy^*}_* \leq \norm{\vx^*}_{X^*} \norm{\vy^*}_{Y^*}, \quad \forall \vx^* \in X^*, ~ \vy^* \in Y^*.
\]
Here, $X^*$ and $Y^*$ denote the norm dual of Banach spaces $X$ and $Y$. 
In the remainder, we will use the following tensor norms. 

The first tensor norm defined by
\[
\|\mT\|_\vee := 
\sup_{\norm{\vx^*}_{X^*} \leq 1, \norm{\vy^*}_{Y^*} \leq 1} |\langle \vx^* \otimes \vy^*, \mT \rangle|
\]
is called the \textit{injective} norm. 
The resulting Banach space equipped with the injective norm is denoted by $X \mathbin{\widecheck{\otimes}} Y$. 
The injectivity implies that if $Z$ is a closed subspace of $X$, then $Z \mathbin{\widecheck{\otimes}} Y$ is a closed subspace of $X \mathbin{\widecheck{\otimes}} Y$. 
This property will play a crucial role in deriving the entropy estimate in Section~\ref{sec:entropy}. 
Furthermore, the injective norm coincides with the operator norm from $X^*$ to $Y$. 

The second tensor norm is the \textit{projective} norm defined by
\[
\|\mT\|_\wedge = \inf \left\{ \sum_{k=1}^n \|\vx_k\|_X \|\vy_k\|_Y \, : \, n \in \mathbb{N}, \, \mT = \sum_{k=1}^n \vx_k \otimes \vy_k \right\}.
\]
The resulting Banach space with the projective norm is denoted by $X \mathbin{\widehat{\otimes}} Y$.
The projectivity implies that if $Z$ is a subspace of $X$, then $(X/Z) \mathbin{\widehat{\otimes}} Y$ is a quotient of $X \mathbin{\widehat{\otimes}} Y$, where $X/Z$ denotes the quotient of $X$ with respect to $Z$. Therefore, there exists a surjection from $X \mathbin{\widehat{\otimes}} Y$ to $(X/Z) \mathbin{\widehat{\otimes}} Y$. 

In a special case where $X = \ell_\infty^K$ and $Y = \ell_2^N$, the injective norm on $\ell_\infty^K \otimes \ell_2^N$ coincides with the max-block-$2$-norm defined by
\[
\|(\vx_1,\dots,\vx_N)\|_{\ell_\infty^K(\ell_2^N)} := \max_{k \in K} \norm{\vx_k}_2, \quad \vx_1,\dots,\vx_K \in \ell_2^N.
\]
The corresponding Banach space is denoted by $\ell_\infty^K(\ell_2^N)$. 
The norm dual of $\ell_\infty^K(\ell_2^N)$, denoted by $\ell_1^K(\ell_2^N)$, is equipped with the norm given by
\[
\|(\vx_1,\dots,\vx_N)\|_{\ell_1^K(\ell_2^N)} := \sum_{k=1}^K \norm{\vx_k}_2, \quad \vx_1,\dots,\vx_K \in \ell_2^N.
\]

\section{Entropy estimate}
\label{sec:entropy}

The main machinery enabling the proof of Theorem~\ref{thm:main} is Maurey's empirical method \cite{carl1985inequalities}, which provides tail bounds on random processes arising in the analysis. 
In this section, we present and prove the key entropy estimate on the linear operators related to the estimator in \eqref{eq:estimator}. 
We first recall the notion of the covering number to state the entropy estimate results. 
For symmetric convex bodies $D$ and $E$, the \emph{covering number} $N(D,E)$ is defined by
\begin{align*}
N(D,E) & := \min \Big\{ l : \exists \vy_1,\dots,\vy_l \in D, \, D \subset \bigcup_{1\le j \le l} (\vy_j + E) \Big\}.
\end{align*}
Then Maurey's empirical method \citep{carl1985inequalities} provides an upper bound on the integral of the square root of the log-covering number for linear operators from $\ell_1^n$. We use a version of this result \cite{junge2020generalized}, summarized as the following lemma. 

\begin{lemma}[{\cite[Lemma~3.4]{junge2020generalized}}]\label{mm}
Let $\mT \in L(\ell_1^n, \ell_\infty^m(\ell_2^d))$. Then
\begin{align*}
& \int_0^{\infty} \sqrt{\ln N(\mT(B_1),\eta B_{\infty,2})} d\eta \\
& \lesssim 
\sqrt{1+\ln (m \vee n)} \, (1+\ln (m \wedge n))^{3/2} \|\mT\|_\vee.    
\end{align*}
\end{lemma}
\noindent Lemma~\ref{mm} considers the case where the range of $\mT$ is $\ell_\infty^m(\ell_2^d))$. 
Note that the upper bound by Lemma~\ref{mm} is independent of the dimension $d$. 
This is a special case of the original result by Carl \cite{carl1985inequalities}, in which the range space is a Banach space of type-$2$. 

We utilize Lemma~\ref{mm} in order to get an entropy estimate with respect to the $\$$-norm. 
The result is obtain in the following two steps. 
The following lemma, proved in Appendix~\ref{sec:proof:lemma:pi2a_nu1}, shows that the $2$-summing norm of $\mX^*$ is equivalent to the projective norm of $\mX$ up to $\sqrt{2N}$. 

\begin{lemma}
\label{lemma:pi2a_nu1}
Let $\mT \in \ell_\infty^K(\ell_2^N) \otimes \ell_2^M$. 
Then $\norm{\mT}_\$$ satisfies
\[
\norm{\mT}_\$ \leq \norm{\mT}_\wedge \leq \sqrt{2N} \, \norm{\mT}_\$.
\]
\end{lemma}

\noindent Lemma~\ref{lemma:pi2a_nu1} implies that the unit $\$$-norm ball is contained in the projective norm ball of radius $\sqrt{2N}$. Then it remains to obtain an upper bound on the entropy integral with respect to the projective norm. 
The result is stated in the following lemma. The proof is provided in Appendix~\ref{sec:proof:lem:entropy_integral}.

\begin{lemma}
\label{lem:entropy_integral}
Let $X = \ell_\infty^K(\ell_2^N) \mathbin{\widehat{\otimes}} \ell_2^M$, $Y = \ell_\infty^m(\ell_2^d)$, and $\mT \in L(X, Y)$. Suppose that $m \leq 2^{NK+M}$. 
Then
\begin{align*}
& \int_0^{\infty} \sqrt{\ln N(\mT(B_X),\eta B_Y)} d\eta \\
& \lesssim \sqrt{1+\ln m \vee (NK+M)} \, (1+\ln m)^{3/2} \|\mT\|_\mathrm{op}.    
\end{align*}
\end{lemma}

\section{Proof of Theorem~\ref{thm:main}}
\label{sec:proof:thm:main}

In this section, we present the proof of Theorem~\ref{thm:main}. By the optimality of $\widehat{\mX}$, we obtain a basic inequality given by
\begin{equation*}
\sum_{l=1}^L \sum_{k=1}^K \left ( y_{l,k} - \langle \mA_{l,k}, \widehat{\mX} \rangle \right )^2 
\leq \sum_{l=1}^L \sum_{k=1}^K \left ( y_{l,k} - \langle \mA_{l,k}, \mX \rangle \right )^2,
\end{equation*} 
which implies 
\begin{equation}
\label{eq:basic_lemma}
\sum_{l=1}^L \sum_{k=1}^K \langle \mA_{l,k}, \widehat{\mX} - \mX \rangle^2 
\leq 
2 \sum_{l=1}^L \sum_{k=1}^K \langle \mA_{l,k}, \widehat{\mX} - \mX \rangle w_{l,k}.
\end{equation}

Recall the constraint set $\kappa(\alpha,\beta)$ is given as the intersection of two norm balls. 
Since $\widehat{\mX} \in \kappa(\alpha,\beta)$, it satisfies $\|\widehat{\mX}\|_{\infty,\mathrm{F}} \leq \alpha$ and $\|\widehat{\mX}\|_\$ \leq \beta$. Furthermore, since $\mX \in \kappa(\alpha,\beta)$, we also have $\|\mX\|_{\infty,\mathrm{F}} \leq \alpha$ and $\|\mX\|_\$ \leq \beta$. 
Since the two norms are sub-additive, we have $\|\widehat{\mX}-\mX\|_{\infty,\mathrm{F}} \leq 2\alpha$ and $\|\widehat{\mX}-\mX\|_\$ \leq 2 \beta$. In other words, we have 
\[
\widehat{\mX}-\mX \in \kappa(2\alpha,2\beta). 
\]

Then a lower-bound (resp. an upper bound) on the left-hand side of \eqref{eq:basic_lemma} (resp. the right-hand side of \eqref{eq:basic_lemma}) is obtained respectively by the following two lemmas, whose proofs are given in Appendix~\ref{sec:proof:lem:sup_quad} and \ref{sec:proof:lemma:ub_cor}. 

\begin{lemma}
\label{lem:sup_quad}
Under the hypothesis of Theorem~\ref{thm:main}, it holds with probability $1-\zeta$ that
\begin{equation}
\label{eq:lem:sup_quad}
\begin{aligned}
& \sup_{\mZ \in \kappa(\alpha,\beta)}
\left| L \norm{\mZ}_\mathrm{F}^2 - \sum_{l=1}^L \sum_{k=1}^K \langle \mA_{l,k}, \mZ \rangle^2 \right| \\
& \lesssim \alpha^2 \ln(2\zeta^{-1}) + \alpha^2 K \left(\rho + \sqrt{\frac{L \ln(2\zeta^{-1})}{K}} \right),
\end{aligned}
\end{equation}
where
\begin{equation}
\label{eq:def_rho}
\rho := \sqrt{\frac{
(\beta/\alpha)^2 N (NK+M) (\ln K)^3}{LK}}.
\end{equation}
\end{lemma}

\begin{lemma}
\label{lemma:ub_cor}
Under the hypothesis of Theorem~\ref{thm:main}, it holds with probability $1-\zeta$ that
\begin{equation}
\label{eq:lemma_ub_cor}
\begin{aligned}
\sup_{\mZ \in \kappa(\alpha, \beta)} \sum_{l=1}^L \sum_{k=1}^K \langle \mA_{l,k}, \mZ \rangle w_{l,k} 
\lesssim 
\alpha \sigma 
\left(
LK \rho + \sqrt{LK\ln(\zeta^{-1})} \right),
\end{aligned}
\end{equation}
where $\rho$ is defined in \eqref{eq:def_rho}.
\end{lemma}

By plugging in the results by these lemmas to \eqref{eq:basic_lemma}, we obtain that \eqref{eq:basic_lemma} implies 
\begin{align*}
\|\widehat{\mX} - \mX\|_\mathrm{F}^2 
& \lesssim 
\frac{\alpha^2 \ln(\zeta^{-1})}{L} 
+ \alpha^2 K \left(\rho + \sqrt{\frac{\ln(\zeta^{-1})}{LK}} \right) \\
& + \alpha K \sigma
\left(
\rho + \sqrt{\frac{\ln(\zeta^{-1})}{LK}} \right).
\end{align*}
Finally, the simplified upper bound in Theorem~\ref{thm:main} is obtained since the first summand in the right-hand side is dominated by the other summands.

\section{Numerical Results}
\label{sec:numerical}

We performed Monte Carlo simulations on synthesized data to study the empirical performance of the tensor-norm-based convex estimator in \eqref{eq:estimator} relative to the spectral method and its refinement via gradient descent \cite{gigi2020spectral}. 
The sensing matrices and measurement noise are generated as in Theorem~\ref{thm:main} so that $\mB_{k,l}$'s are independent copies of a random matrices whose entries are drawn i.i.d. from $\mathcal{N}(0,1)$ and $w_{k,l}$'s are i.i.d. $\mathcal{N}(0,\sigma^2)$. 
The ground-truth matrix is generated as a rank-$r$ matrix given by $\mU \in \R^{M \times r}$ uniformly distributed on a Stiefel manifold and $\mV_k$'s are independent copies of a random matrix with i.i.d. standard Gaussian entries. 
The convex estimator uses the estimates of the parameters $\alpha$ and $\beta$ given by the corresponding norms computed from the rank-$r$ approximation of 
\begin{equation*}
\widehat{\mX}_0 = \frac{1}{L}\sum_{i=1}^L \sum_{k=1}^K y_{l,k} \left( \bm{e}_k^* \otimes \mB_{k,l} \right) \in \mathbb{R}^{M\times NK}.
\end{equation*}
Convex programs for both the $\$$-norm computation and the convex estimator are implemented as ADMM algorithms, which are derived in Appendix~\ref{sec:admm}. We observe the median estimation error from $20$ instances in the Monte Carlo simulations.

We first compare the estimates of the ground-truth column space respectively by the convex estimator and the spectral method \cite{gigi2020spectral}. 
The error is measured by the sine of the largest principal angle between two subspaces. 

   \begin{figure}[!htb]
           \centering
           \subfigure[Convex estimator]{
   		\includegraphics[scale=0.25]{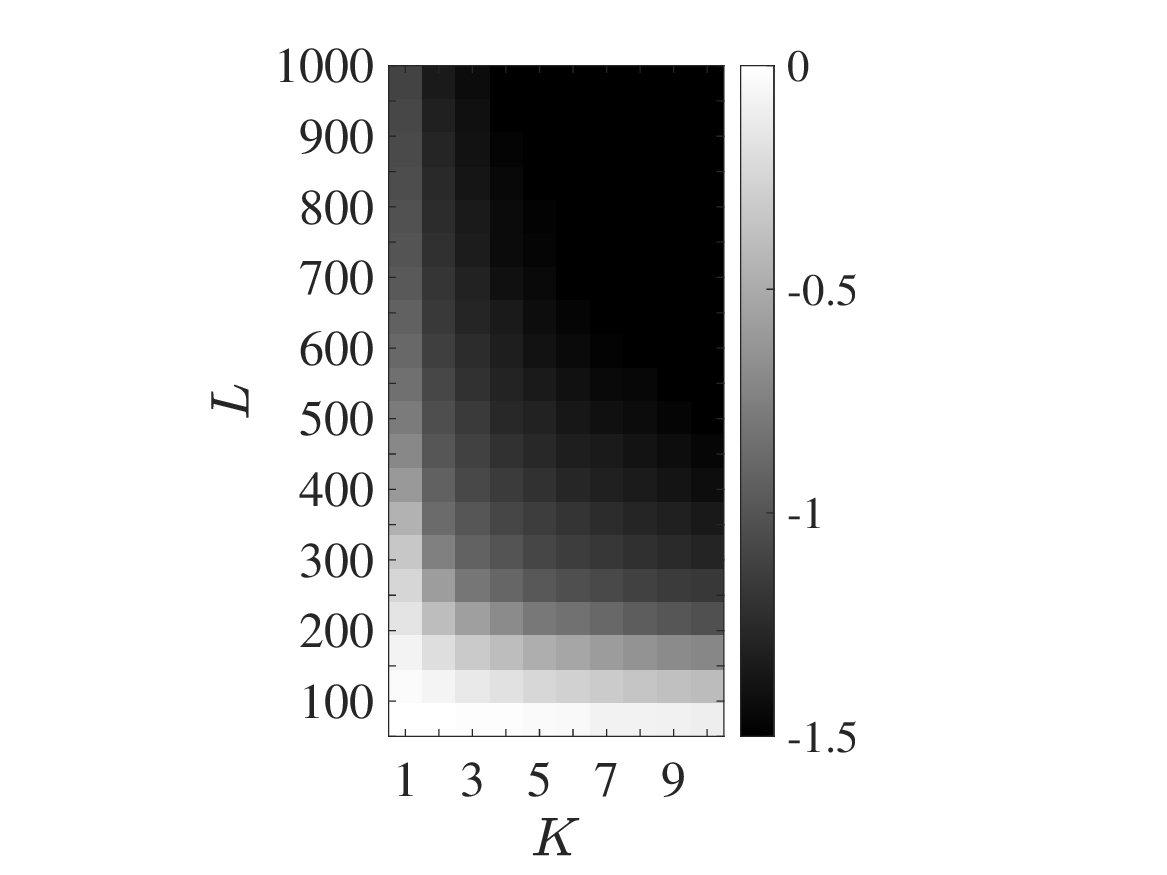}}
   	    \hspace{-0.1\textwidth}
           \subfigure[Spectral estimator]{
   	    \includegraphics[scale=0.25]{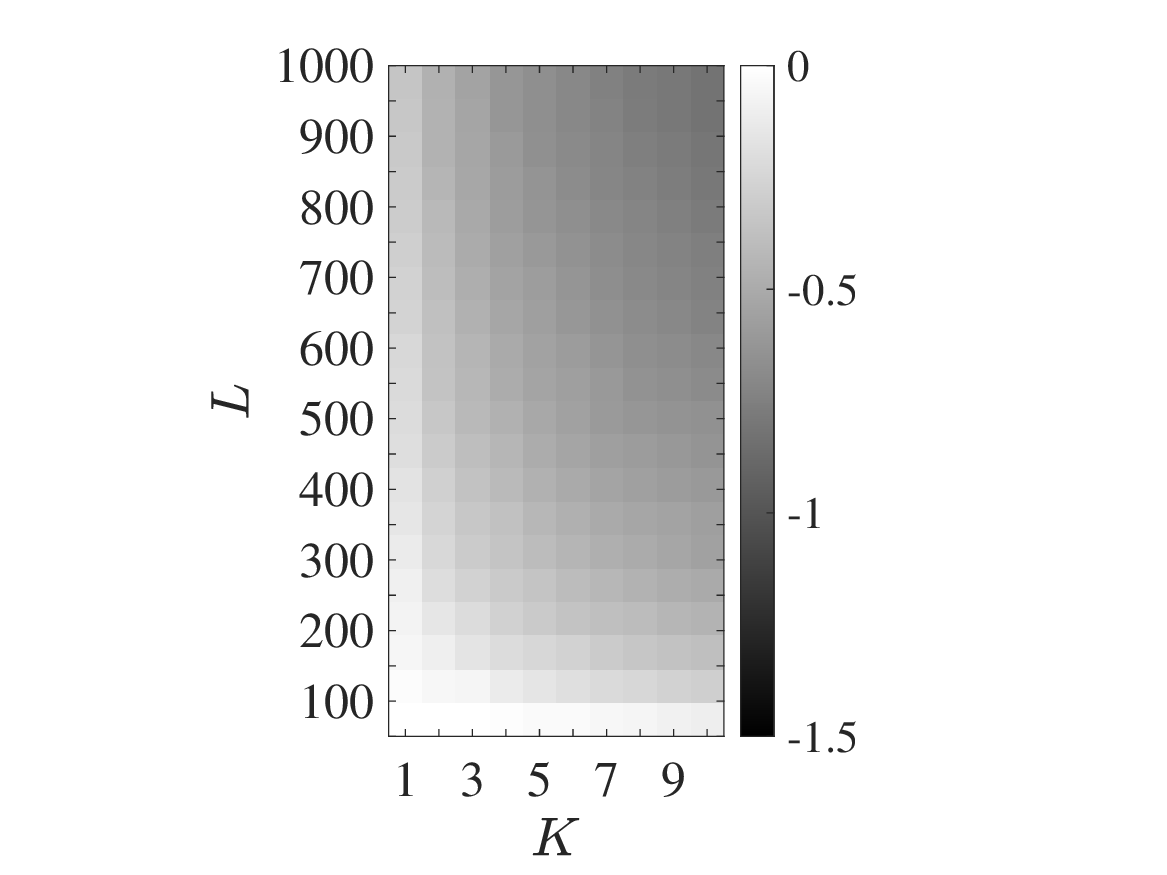}}
           \caption{The log-base-$10$ of the estimation error of the ground-truth column space in the noisy case  ($\mathrm{SNR}=20$dB, $M=100,\,N=20,\,r=2$).  }
           \label{fig:column space_noise} 
   \end{figure}

\begin{figure}[!htb]
           \centering
           \subfigure[Convex estimator]{
   		\includegraphics[scale=0.25]{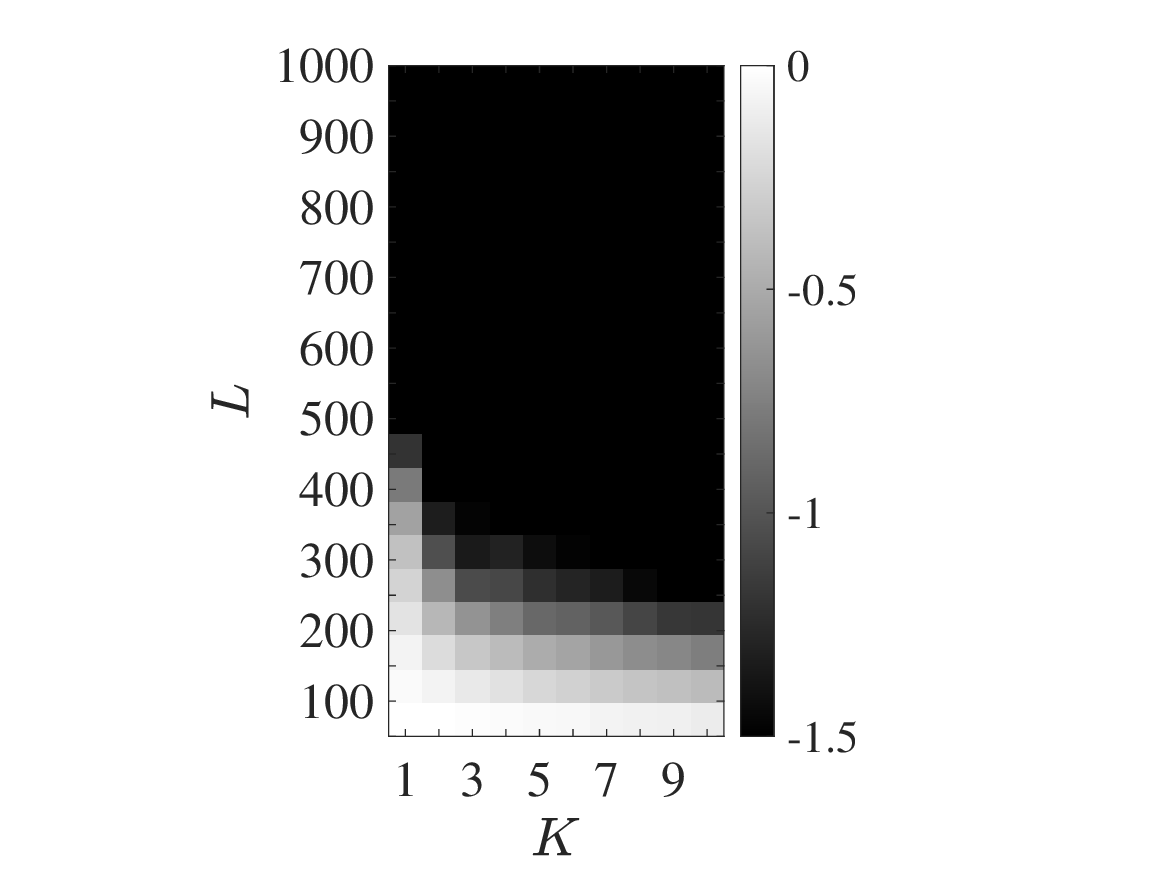}}
   	    \hspace{-0.1\textwidth}
           \subfigure[Spectral estimator]{
   	    \includegraphics[scale=0.25]{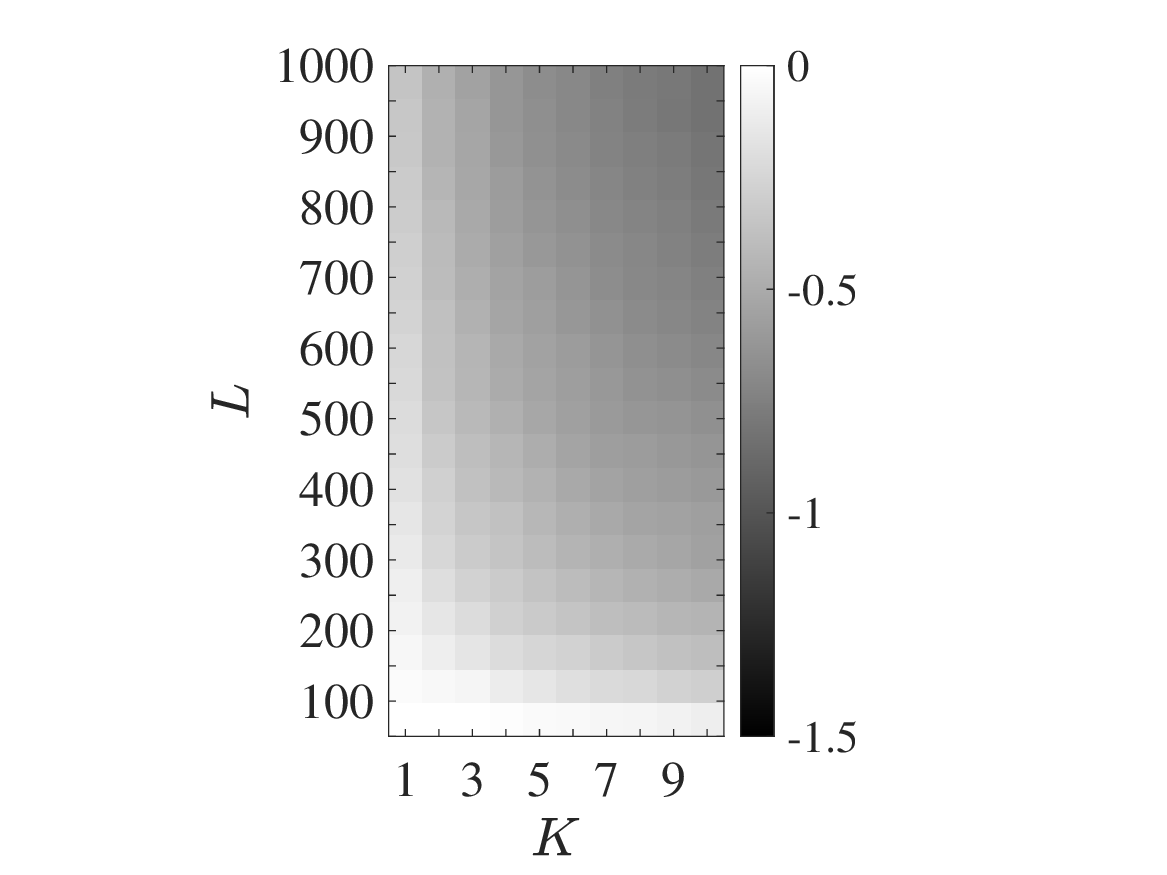}}
           \caption{The log-base-$10$ of the estimation error of the ground-truth column space in the noiseless case ($M=100,\,N=20,\,r=2$).  }
           \label{fig:column space_noiseless} 
   \end{figure}

\Cref{fig:column space_noise} compares the estimation error by the convex estimator and the spectral method in the noisy case with SNR $20$dB.
The errors by both estimators decay with larger $L$ and $K$. 
However, in all observed regime of the parameters, the convex estimator outperforms the spectral estimator. 
As shown in \Cref{fig:column space_noiseless}, the comparison between the two estimator remains similar in the noiseless case. 

Next we compare the performance of estimating the entire ground-truth matrix $\bm{X}$ by the convex estimator and the gradient descent from spectral initialization \cite{gigi2020spectral}. 
In this comparison, the metric is chosen as the normalized reconstruction error given by ${\|\widehat{\mX}-\mX\|_{\mathrm{F}}^2}/{\|\mX\|_{\mathrm{F}}^2}$, where $\widehat{\mX}$ denotes an estimate of $\bm{X}$. 
  
   \begin{figure}[!htb]
           \centering
           \subfigure[Convex estimator]{
   		\includegraphics[scale=0.25]{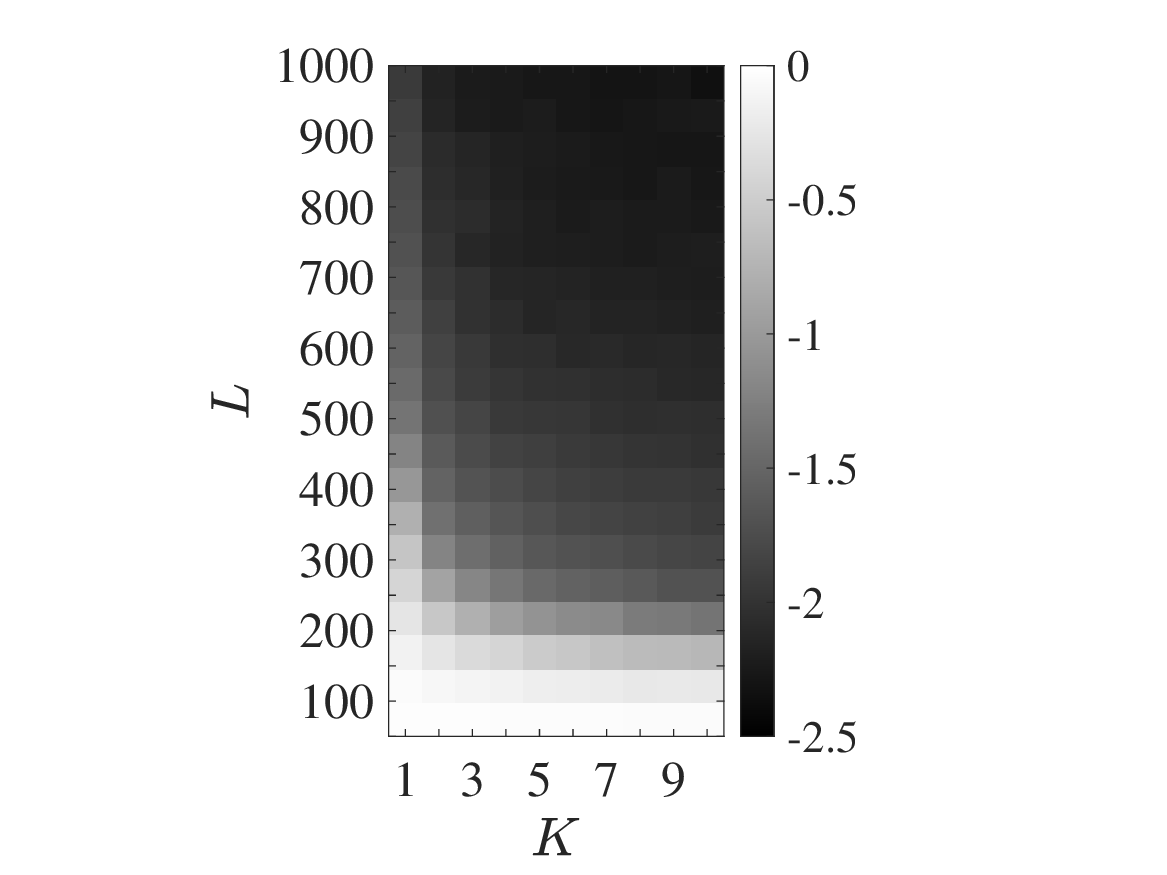}
   	    \hspace{-0.1\textwidth}}
           \subfigure[Gradient descent]{
   		\includegraphics[scale=0.25]{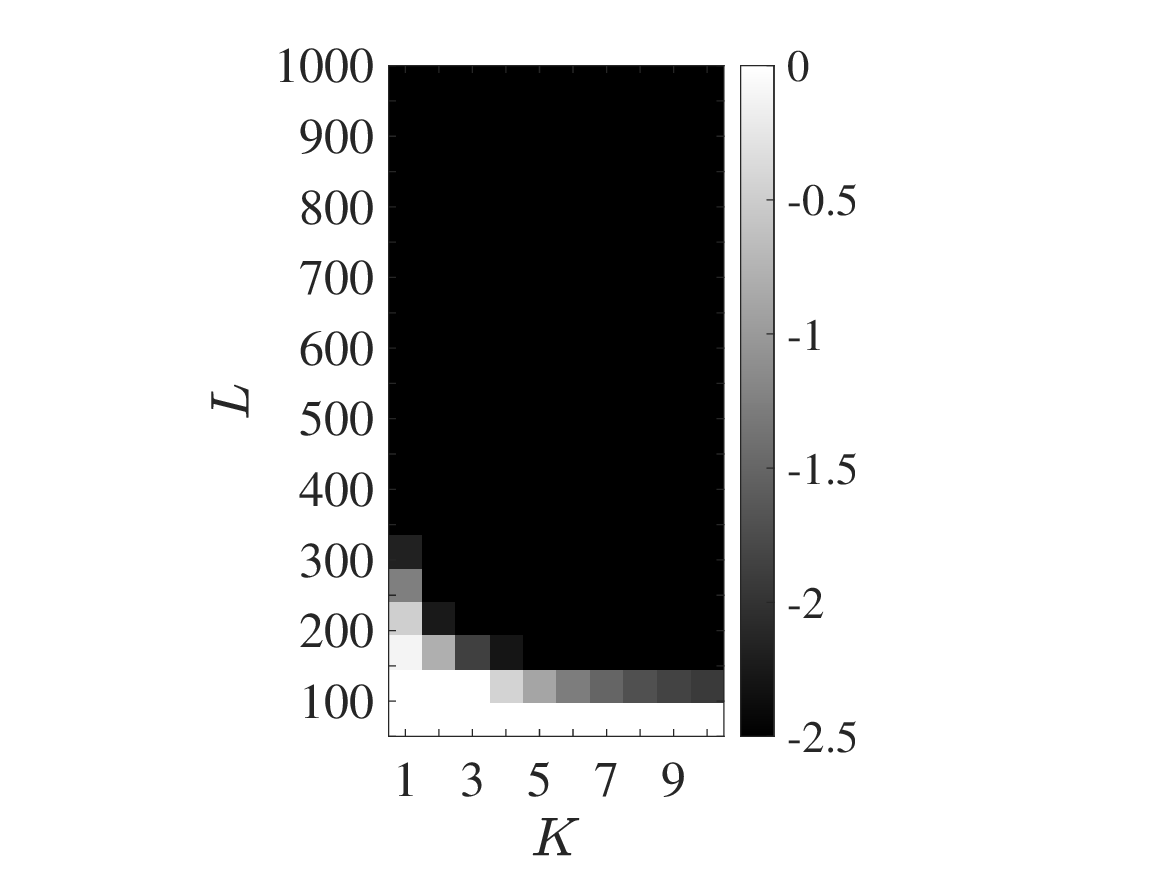}}
           \caption{The log-base-$10$ of the estimation error of the ground-truth matrix ($\mathrm{SNR}=20$ dB, $M=100,\,N=20,\,r=2$).}
           \label{fig:phase_noise} 
  \end{figure}

\begin{figure}[!htb]
           \centering
           \subfigure[Convex estimator]{
   		\includegraphics[scale=0.25]{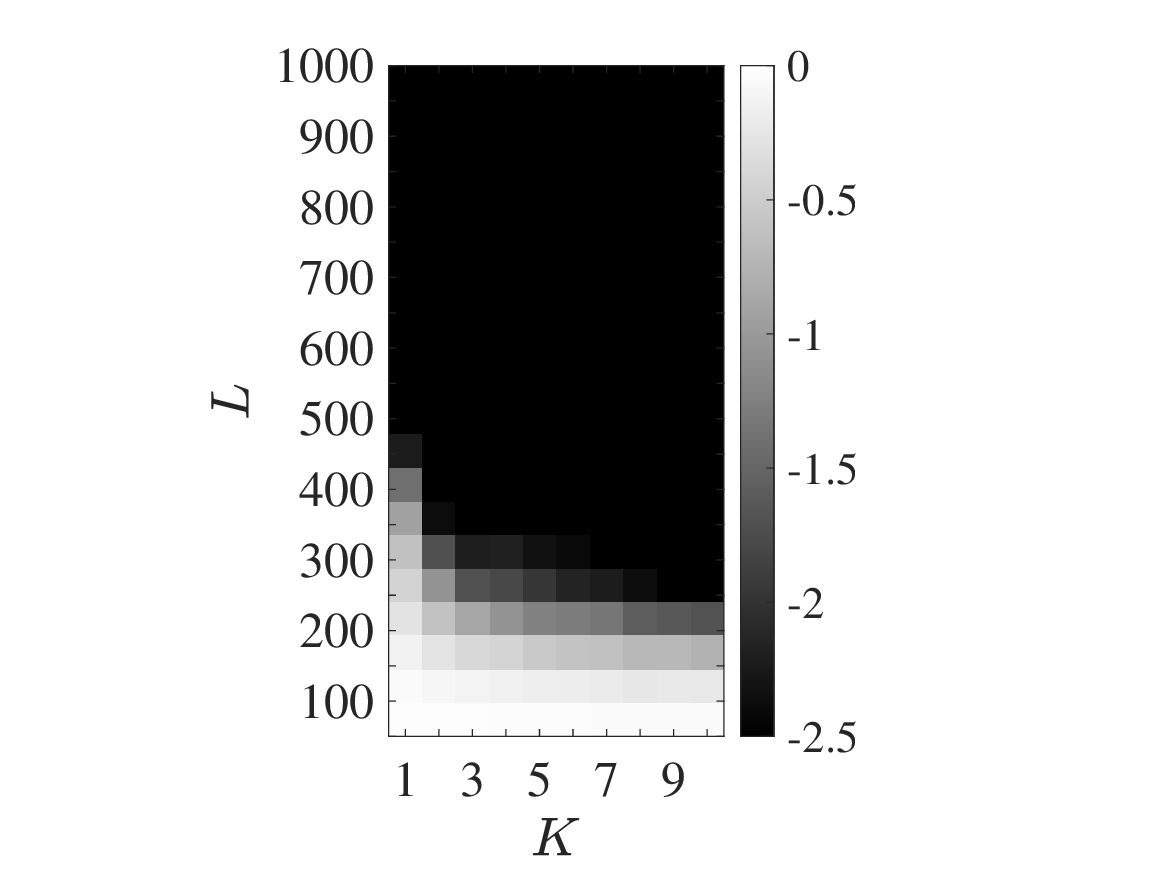}
   	    \hspace{-0.1\textwidth}}
           \subfigure[Gradient descent]{
   		\includegraphics[scale=0.25]{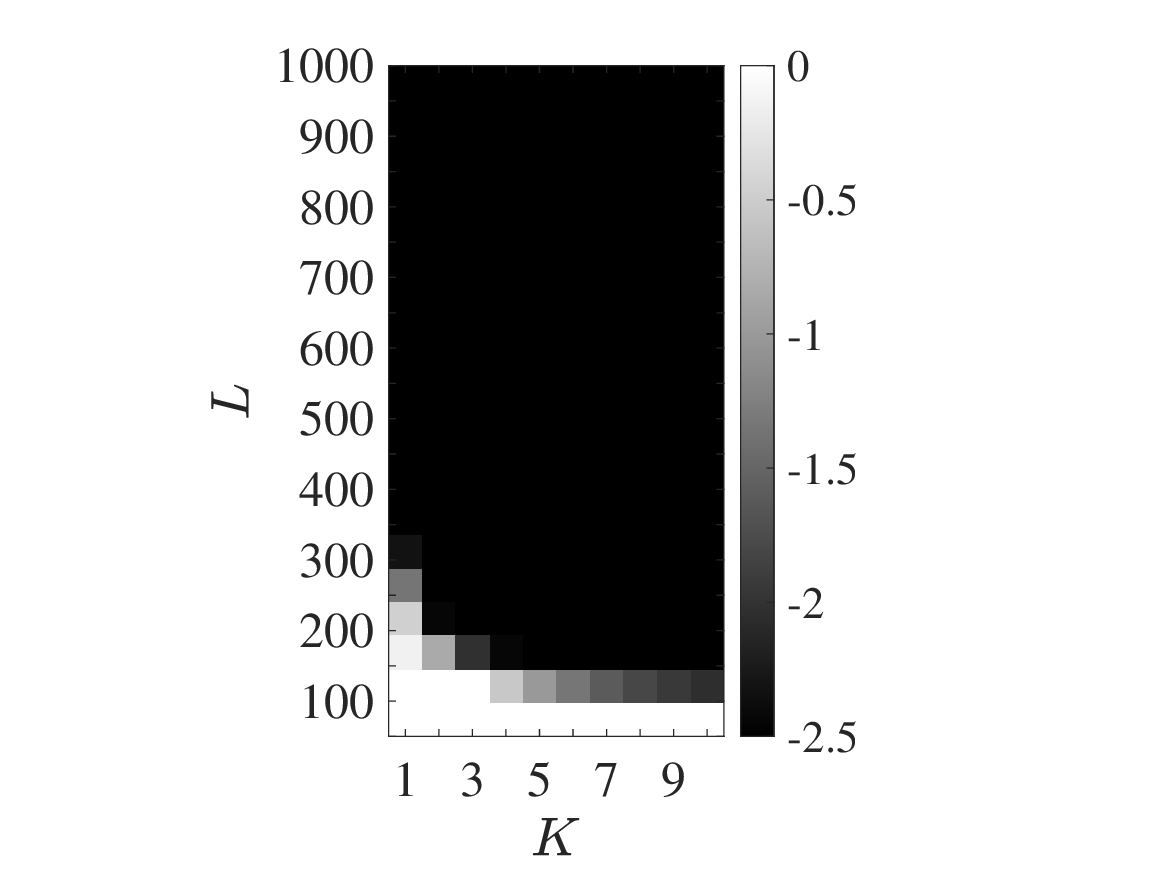}}
           \caption{The log-base-$10$ of the estimation error of the ground-truth matrix in the noiseless case ($M=100,\,N=20,\,r=2$).}
           \label{fig:phase_noiseless} 
  \end{figure}

\Cref{fig:phase_noise} demonstrates that the gradient descent provides a better empirical phase transition than the convex estimator. 
However, while our main result provides a rigorous estimation error bound for the convex estimator, such a theoretical analysis of the gradient descent method has yet to be established.
For both estimators, the error decays with larger $K$ and $L$. 
The phase transition between success (error $\leq 10^{-1.5}$) and failure by the convex estimator occurs on a boundary in which the threshold on $L$ decays with $K$ until $M/K$ is dominated by $N$. 
This corroborates the theoretical analysis in Theorem~\ref{thm:main}. 
Furthermore, unlike the result in Theorem~\ref{thm:main}, the estimation error by the convex estimator continues to decrease with higher SNR. 
As shown in \Cref{fig:phase_noiseless}, the normalized estimation error is below $10^{-2.5}$ when $L$ is above the displayed threshold. 
The convex estimator provides a significantly improved estimation performance in the noiseless case. In particular, the phase transition by the convex estimator is comparable to that by the gradient descent method. 

\section{Conclusion}
\label{sec:conclusion}
In this paper, we proposed a convex program that jointly estimates a set of low-rank matrices sharing a common column space. 
The data model arises in many applications in practice including remote sensing, multi-class learning, and multi-spectrum imaging. 
The estimation problem is equivalently rewritten as block-wise sensing of a low-rank matrix. 
We have shown that the proposed convex estimator leads to a more favorable sample complexity than the individual recovery of each block when the number of blocks $K$ and the dimension $M$ are high relative to the number of columns per block $N$. 
We provide an ADMM algorithm to tackle large-sized problems. 
In the future work, we will investigate the performance of the convex estimator with faster sketching models via fast Johnson-Lindenstrauss transforms \cite{woodruff2014sketching}.


\appendix

\subsection{Concentration Inequalities}
\label{sec:concentration}

We use a set of concentration inequalities for the proofs in this paper. The first lemma provides a tail bound on the $\ell_1$-norm of an image of an arbitrary column vector via a Gaussian random matrix. 

\begin{lemma}[{\cite[Lemma~2.1]{plan2014dimension}}]
\label{lem:ell1embed}
Let $\va_1, \dots, \va_m$ be independent copies of $\va \sim \mathcal{N}(\vzero, \mI_n)$. 
Let $K \subset \R^n$ be a bounded subset. Then
\begin{align*}
& \left| \frac{1}{m} \sum_{i=1}^m \left| \langle \va_i, \vx \rangle \right| - \sqrt{\frac{2}{\pi}} \norm{\vx}_2 \right| \\
& \leq \frac{4w(K)}{\sqrt{m}} + \frac{d(K)\sqrt{2\ln(2\zeta^{-1})}}{\sqrt{m}}    
\end{align*}
holds with probability at least $1-\zeta$, where $w(K)$ denote the Gaussian width of $K$ and $d(K) = \sup_{\vx \in K} \norm{\vx}_2$. 
\end{lemma}

\begin{remark} 
Lemma~\ref{lem:ell1embed} implies that there exists a numerical constant $c$ such that $\ell_2^n$ is embedded into a subspace of $\ell_1^m$ via some $\mPhi\in \R^{m \times n}$ with $m = c \delta^{-2} n$ so that $|\norm{\vx}_2 - \norm{\mPhi \vx}_1| \leq \delta$ for all $\vx \in \mathbb{S}^{n-1}$.
\end{remark}

\noindent The next lemma is a consequence of Dudley's inequality and provides a tail bound on the supremum of a Gaussian random process. 
\begin{lemma}[{\cite[Theorem~8.1.6]{vershynin2018high}}]
\label{lemma:gaussianmax}
Let $\vxi \sim \mathcal{N}(\vzero,\mI_n)$, $\Delta \subset \mathbb{R}^n$, and $\zeta \in (0,1)$. 
Then 
\[ \sup_{\vf \in \Delta}  \left | \vf^* \vxi \right | 
\lesssim \int_{0}^{\infty} \sqrt{\ln N(\Delta, \eta B_2)} d\eta + \mathrm{diam}(\Delta) \sqrt{\ln (\zeta^{-1}) }
\] 
holds with probability $1 - \zeta$, where $\mathrm{diam}(\Delta)$ denotes the diameter of $\Delta$ in $\ell_2^n$. 
\end{lemma}

\noindent We also use the results on the suprema of second-order chaos processes \cite{krahmer2014suprema}, summarized as the following theorem. 

\begin{theorem}[{Theorem~3.1 in \citep{krahmer2014suprema}}]
\label{thm:kmr}
Let $\vxi \in \mathbb{R}^n$ be a Gaussian vector with $\mathbb{E}[\vxi] = \vzero$ and $\mathbb{E}[\vxi \vxi^*] = \mI_n$. 
Let $\Delta \subset \mathbb{R}^{m \times n}$. 
Then 
\begin{align*}
& \sup_{\mQ \in \Delta} \left|\norm{\mQ \vxi}_2^2 - \E[\norm{\mQ \vxi}_2^2] \right| \\
& \lesssim E + V \sqrt{\ln(2\zeta^{-1})} + U \ln(2\zeta^{-1})
\end{align*}
holds with probability $1-\zeta$, where
\begin{align*}
E {} & := \gamma_2(\Delta)
\left[ \gamma_2(\Delta) + d_{\mathrm{F}}(\Delta) \right], \\
V {} & := d_{\mathrm{S}}(\Delta) \left[ \tilde{\gamma}_2(\Delta) + d_{\mathrm{F}}(\Delta) \right], \\
U {} & := d_{\mathrm{S}}^2(\Delta).
\end{align*}
Here $\gamma_2(\Delta)$ denotes the Talagrand $\gamma_2$-functional of the metric space given by the spectral norm, and $d_\mathrm{S}(\Delta)$ and $d_\mathrm{F}(\Delta)$ denotes the radii of $\Delta$ with respect to the spectral norm and the Frobenius norm, respectively.
\end{theorem}

\subsection{Embedding $\ell_1^N$ to $\ell_\infty^{2^N}$}
\label{sec:ell1_to_ellinf}

Let $\iota: \ell_1^N \to \ell_\infty^{2^N}$ denote a linear map defined by
\[
\iota\left((x_n)_{n=1}^N\right) = \left(\sum_{n=1}^N \epsilon_n x_n\right)_{(\epsilon_n)_{n=1}^N \in \{\pm 1\}^N}.
\]
Then $\ell_1^N$ is isometrically embedded into $\ell_\infty^{2^N}$, i.e. 
\begin{align*}
\left\|\iota((x_n)_{n=1}^N)\right\|_\infty
& = \max \left\{ \bigg| \sum_{n=1}^N \epsilon_n x_n \bigg| : (\epsilon_n)_{n=1}^N \in \{\pm 1\}^N \right\} \\
& = \left\|(x_n)_{n=1}^N\right\|_1.
\end{align*}

Let $X = \ell_\infty^{2^N}$ and $E = \iota(\ell_1^N) \subset X$. 
Then the dual space of linear functionals on $X$ is denoted by $X^* = \ell_1^{2^N}$. 
We denote the vector space of linear functionals on $E$ by $E^*$. 
First we note that the restriction of $\iota^*$ on $E^*$ is an isometric bijection. Indeed, we have
\begin{align*}
\norm{\vy}_{E^*} 
& = \sup_{\norm{\iota(\vx)}_\infty \leq 1} \langle \iota(\vx), \vy \rangle
= \sup_{\norm{\vx}_1 \leq 1} \langle \vx, \iota^*(\vy) \rangle \\
& = \norm{\iota^*(\vy)}_\infty, \quad \forall \vy \in E^*.     
\end{align*}
Next, due to the Hahn-Banach theorem, for any $\vy \in E^*$, there exists a linear functional $\tilde{\vy} \in X^*$ such that $\|\tilde{\vy}\|_{X^*} = \|\vy\|_{E^*}$. 
Consequently, there exists an isometric bijection map $\varrho$ from $X^*/E^\perp$ to $E^*$, where $E^\perp = \{ \vy \in X^* : \langle \vy, \vx \rangle = 0, \forall \vx \in E\}$ and $\langle \cdot, \cdot \rangle$ denotes the canonical bilinear transform on $X^* \times X$. 
Therefore, the map $\iota^* \circ \varrho: X^*/E^\perp \to \ell_\infty^N$ is an isometric bijection. 
Furthermore, the quotient map $q: X^* \to X^*/E^\perp$ is a metric surjection \cite[Eqs. (1.3.2) and (2.2.5)]{carl2008entropy}, i.e. 
\[
q(B_{X^*}) = B_{X^*/E^\perp},
\]
where $B_{X^*}$ and $B_{X^*/E^\perp}$ denote the unit norm ball respectively in $X^*$ and $X^*/E^\perp$. 
Finally, we deduce that $\iota^* \circ \varrho \circ q: \ell_1^{2^N} \twoheadrightarrow \ell_\infty^N$ is a metric surjection.

\subsection{Duality} 
\label{sec:duality}

Let $X$ and $Y$ be finite-dimensional Banach spaces. Let $\mA: X \to Y$ be a linear operator such that 
\begin{equation}
\label{eq:embed_A}
\left| \norm{\mA \vx}_Y - \norm{\vx}_X \right| \leq \delta \norm{\vx}_X, \quad \forall \vx \in X. 
\end{equation}
Let $E = \mA(X)$ denote the image of $X$ via $\mA$, which is a subspace of $Y$, i.e. $X \hookrightarrow E \subset Y$. 
Let $E^\perp := \{\vy \in Y^* : \langle \vy, \vx \rangle = 0, \forall \vx \in E\}$. 
Then, by the Hahn-Banach theorem, $E^*$ is isometrically isomorphic to $Y^*/E^\perp$ and there exists an isometric bijection $\varrho$ from $Y^*/E^\perp$ to $E^*$. 
It follows from \eqref{eq:embed_A} that $\mB = \mA^*|_{E^*}$ satisfies
\begin{equation*}
\left| \norm{\mB \vy}_{X^*} - \norm{\vy}_{E^*} \right| \leq \delta \norm{\vy}_{E^*}, \quad \forall \vy \in E^*. 
\end{equation*}
Therefore, we obtain that $\mB \circ \varrho: Y^*/E^\perp \twoheadrightarrow X^*$ is a bijection satisfying $\norm{\mB \circ \varrho} \leq 1+\delta$. 
Furthermore, the quotient map $q: Y^* \to Y^*/E^\perp$ is a metric surjection.

\subsection{Proof of Lemma~\ref{lem:interlacing}}
\label{sec:proof:lem:interlacing}

Let $\mU \in \R^{M \times r}$ and $\mV \in \R^{MK \times r}$ satisfy that $\mX = \mU \mV^*$ and $\mU^* \mU = \mI_r$. Then we have
\begin{align*}
\norm{\mX}_\$ 
& \leq \norm{\mU}_\mathrm{F} \norm{\mV^*}_{\infty,\mathrm{F}}
\leq \sqrt{r} \norm{\mU} \norm{\mV^*}_{\infty,\mathrm{F}} \\
& = \sqrt{r} \norm{\mU \mV^*}_{\infty,\mathrm{F}}
= \sqrt{r} \norm{\mX}_{\infty,\mathrm{F}},
\end{align*}
which implies the upper bound in \eqref{eq:interlacing}.

To derive the lower bound in \eqref{eq:interlacing}, we consider $\mU$ and $\mV$ satisfy that $\mX = \mU \mV^*$ and $\norm{\mX}_\$ = \norm{\mU}_\mathrm{F} \norm{\mV}_{\infty,\mathrm{F}}$. 
The common number of columns of $\mU$ and $\mV$ is not necessarily $r$ this time. 
Let $\mX_k = \mX (\ve_k \otimes \mI_N)$ and $\mV_k^* = \mV^* (\ve_k \otimes \mI_N)$ denote the $k$th block of $\mX$ and $\mV^*$, respectively. 
Let $k_\star = \mathop{\mathrm{argmax}}_{k \in [K]} \norm{\mX_k}_\mathrm{F}$. 
Then
\begin{align*}
\norm{\mX}_{\infty,\mathrm{F}} 
& = \norm{\mX_{k_\star}}_\mathrm{F} 
= \norm{\mU \mV_{k_\star}^*}_\mathrm{F} 
\leq \norm{\mU}_\mathrm{F} \norm{\mV_{k_\star}^*}_\mathrm{F} \\
& = \norm{\mU}_\mathrm{F} \norm{\mV^*}_{\infty,\mathrm{F}} 
= \norm{\mX}_\$. 
\end{align*}
This completes the proof. 

\subsection{Proof of Lemma~\ref{lemma:pi2a_nu1}}
\label{sec:proof:lemma:pi2a_nu1}

We first show that the $\$$-norm is a valid tensor norm. 
Let $x \in \ell_\infty^K(\ell_2^N)$ and $y \in \ell_2^M$. 
Since $x \otimes y$ is rank-$1$, the optimal factorization in the definition of the $\$$-norm is through the trivial $1$-dimensional space and hence it follows that
\[
\norm{x \otimes y}_\$ = \norm{x}_{\ell_\infty^K(\ell_2^N)} \norm{y}_{\ell_2^M}.
\]
Similarly, for $x^* \in \ell_1^K(\ell_2^N)$ and $y^* \in \ell_2^M$, the dual norm of $\$$-norm on $x^* \otimes y^*$ is written as
\begin{align*}
\sup_{\begin{subarray}{c} \norm{x}_{\ell_\infty^K(\ell_2^N)} \leq 1 \\  \norm{y}_{\ell_2^M} \leq 1 \end{subarray}}
\langle (x^* \otimes y^*) x, y \rangle
& = \sup_{\begin{subarray}{c} \norm{x}_{\ell_\infty^K(\ell_2^N)} \leq 1 \\  \norm{y}_{\ell_2^M} \leq 1 \end{subarray}}
\langle x^*, x \rangle \, \langle y^*, y \rangle \\
& = \norm{x^*}_{\ell_1^K(\ell_2^N)} \norm{y^*}_{\ell_2^M},
\end{align*}
where $\langle x^*, x \rangle$ denotes the dual bracket representing the evaluation of the linear functional $x^*$ on $x$. 
Therefore, we have shown that the $\$$-norm is a tensor norm. 
Then, since the projective norm is the largest tensor norm, it follows that $\norm{\mT}_\$ \leq \norm{\mT}_\wedge$. 

To show the remaining inequality $\norm{\mT}_\wedge \leq \sqrt{2N} \norm{\mT}_\$$, we use the $2$-summing norm of of $\mT \in X \otimes Y$ defined as the smallest constant $c > 0$ that satisfies
\[
\sum_k \|\mT \vx^*_k\|_Y^2
\leq c^2 \,  \sup_{\norm{\vx}_{X^{**}} \leq 1} \sum_k |\langle \vx, \vx^*_k\rangle|^2
\]
for all sequences $(\vx^*_k) \subset X^*$. Here, $X^{**}$ denotes the double dual, which coincides with $X$ since we consider the finite-dimensional case. 
The $2$-summing norm will be denoted by $\pi_2(\mT)$. 
In the finite-dimensional case, the $2$-summing norm $\pi_2$ is self-dual by satisfying 
\begin{equation}
\label{eq:pi2_self_dual}
\mathrm{tr}(\mS \mT^*) \leq \pi_2(\mS) \pi_2(\mT^*).
\end{equation}
To be self-contained, below we present the derivation of the inequality in \eqref{eq:pi2_self_dual}. 
The arguments are taken from \cite{jameson1987summing}. 
We first recall $\mu_p$ defined on a finite sequence $(x_1,\dots,x_k)$ in a normed space $X$ as 
\[
\mu_p(x_1,\dots,x_k) := \sup \left\{ \left(\sum_{j=1}^k |f(x_j)|^p\right)^{1/p} : f \in B_{X^*} \right\},
\]
where $B_{X^*}$ denotes the unit ball in the dual space $X^*$. 
Then $p$-nuclear norm of a linear operator $\mT$ from a normed space $X$ to another normed space $Y$ is defined by
\begin{align*}
& \nu_p(\mT) \\
& := \inf \left\{ \left(\sum_{i=1}^k \|f_i\|^p\right)^{1/p} \mu_{p'}(y_1,\dots,y_k) : \mT = \sum_{i=1}^k f_i \otimes y_i \right\},    
\end{align*}
where $p'$ satisfies $1/p + 1/p' = 1$. 
Then by the definition of the trace, we have
\begin{equation}
\label{eq:tr_nu1}
\mathrm{tr}(\mS \mT^*) \leq \nu_1(\mS \mT^*).
\end{equation}
Moreover, by the definition of the nuclear norm and $2$-summing norm, it has been shown \cite[4.2]{jameson1987summing} that
\begin{equation}
\label{eq:nu1_pi2nu2}
\nu_1(\mS \mT^*) \leq \pi_2(\mS) \nu_2(\mT^*).  
\end{equation}
Finally, since we consider the finite-dimensional case, the $2$-summing norm and the $2$-nuclear norm coincides \cite[Theorem 5.11]{jameson1987summing}. 
Therefore, the inequality in \eqref{eq:pi2_self_dual} follows from \eqref{eq:tr_nu1} and \eqref{eq:nu1_pi2nu2}. 

Armed with the inequality in \eqref{eq:pi2_self_dual}, we proceed to the remainder of the proof of Lemma~\ref{lemma:pi2a_nu1}. 
By the trace duality, the projective norm of $\mT^* \in \ell_2^M \otimes \ell_\infty^K(\ell_2^N) = L(\ell_2^M,\ell_\infty^K(\ell_2^N))$ is written as
\begin{equation}
\label{eq:pi2_dual_blk}
\norm{\mT^*}_\wedge = \sup \{ \mathrm{tr}(\mS \mT^*) : \mS \in L(\ell_\infty^K(\ell_2^N),\ell_2^M), \, \|\mS\|_\vee \leq 1 \}.
\end{equation}
Then, by the trivial decomposition of $\mS = \mS \circ \mathrm{id}$ via $\ell_\infty^{NK}$, we have
\begin{equation}
\label{eq:pi2_S*1}
\begin{aligned}
& \pi_2(\mS: \ell_\infty^K(\ell_2^N) \to \ell_2^M) \\
& \leq \| \mathrm{id}: \ell_\infty^K(\ell_2^N) \to \ell_\infty^{NK} \| \cdot \pi_2(\mS: \ell_\infty^{NK} \to \ell_2^M) \\
& \leq \sqrt{2} \, \|\mS: \ell_\infty^{NK} \to \ell_2^M\| \\
& \leq \sqrt{2} \, \| \mathrm{id}: \ell_\infty^{NK} \to \ell_\infty^K(\ell_2^N) \| \cdot \|\mS: \ell_\infty^K(\ell_2^N) \to \ell_2^M\| \\
& \leq \sqrt{2N} \, \|\mS: \ell_\infty^K(\ell_2^N) \to \ell_2^M\|,
\end{aligned}
\end{equation}
where the second inequality follows from \cite[Propositions~9.3 and 9.8]{jameson1987summing}. 
By plugging in \eqref{eq:pi2_self_dual} and \eqref{eq:pi2_S*1} into \eqref{eq:pi2_dual_blk}, we obtain
\[
\norm{\mT^*}_\wedge \leq \sqrt{2N} \, \pi_2(\mT^*).
\]
Furthermore, since all Banach spaces here are finite-dimensional, it follows from \cite[Proposition~1.13]{jameson1987summing} that
$\norm{\mT}_\wedge = \norm{\mT^*}_\wedge$. 
Therefore, we have shown that
\begin{equation}
\label{eq:proj_ub_pi2}
\norm{\mT}_\wedge \leq \sqrt{2N} \, \pi_2(\mT^*).
\end{equation}

The following lemma provides an alternative characterization of the $2$-summing norm so that one can compare the $\$$-norm and the $2$-summing norm on the dual of $\ell_\infty^K(\ell_2^N) \otimes \ell_2^M$.

\begin{lemma}[{\cite[Lemma~3.3]{lee2021approximately}}]
\label{lemma:p2ANDp}
Let $\mT \in L(X^*,Y)$ with $X$ complete. Then the $2$-summing norm of the adjoint $\mT^*$ is expressed as
\begin{align*}
\pi_2(\mT^*)
:= \inf & \{ \pi_2(\mT_1^*) \|\mT_2^*\|_\vee \,:\, d \in \mathbb{N}, \, \mT_1^* \in L(Y^*,\ell_2^d), \, \\ 
& ~ \mT_2^* \in L(\ell_2^d,X), \, \mT^* = \mT_2^* \mT_1^* \}.
\end{align*}
\end{lemma}

\noindent Let $\mT \in \ell_\infty^K(\ell_2^N) \otimes \ell_2^M = L(\ell_1^K(\ell_2^N),\ell_2^M)$ be factorized as $\mT = \mU \mV^*$ via $\ell_2^d$ with $\mV^* \in L(\ell_1^K(\ell_2^N),\ell_2^d)$ and $\mU \in L(\ell_2^d,\ell_2^M)$ for some $d \in \mathbb{N}$.  
Then it follows that 
\[
\norm{\mV}_\vee 
= \norm{\mV^*}_\vee
= \norm{\mV^*}_{\infty,\mathrm{S}}, 
\]
where $\norm{\cdot}_{\infty,\mathrm{S}}$ denotes the maximum-block-spectral norm defined by
\begin{align*}
\norm{[\mV_1^* ~ \mV_2^* ~ \cdots ~ \mV_K^*]}_{\infty,\mathrm{S}}
& = \norm{\sum_{k=1}^K \ve_k^* \otimes \mV_k^*}_{\infty,\mathrm{S}} \\
& = \max_{k \in [K]} \norm{\mV_k^*}. 
\end{align*}
Furthermore, since $\mU^* \in L(\ell_2^d,\ell_2^M) = \ell_2^d \otimes \ell_2^M$, the $2$-summing norm and Frobenius norm of $\,U^*$ coincide, i.e. 
\[
\pi_2(\mU^*) = \norm{\mU^*}_\mathrm{F} = \norm{\mU}_\mathrm{F}.
\]
Therefore, the $2$-summing norm of $\mT^*$ is written as
\[
\pi_2(\mT^*) = \inf_{\mU,\mV: \mU\mV^* = \mT} \norm{\mU}_\mathrm{F} \|\mV^*\|_{\infty,\mathrm{S}},
\]
Then we deduce that the $\$$-norm defined in in \eqref{eq:dollar_norm} satisfies
\[
\pi_2(\mT^*) \leq \norm{\mT}_\$,
\]
which together with \eqref{eq:proj_ub_pi2} implies $\norm{\mT}_\wedge \leq \sqrt{2N} \norm{\mT}_\$$. 
This completes the proof.

\subsection{Proof of Lemma~\ref{lem:entropy_integral}}
\label{sec:proof:lem:entropy_integral}
Recall that the dual of $\ell_\infty^K(\ell_2^N)$ is given by $(\ell_\infty^K(\ell_2^N))^* = \ell_1^K(\ell_2^N)$. 
Similarly, the range space $\ell_\infty^m \mathbin{\widecheck{\otimes}} \ell_2^d$ is identified to $\ell_\infty^m(\ell_2^d)$. 
Due to Lemma~\ref{lem:ell1embed} in Appendix~\ref{sec:concentration} (also see the remark after the lemma), there exists a linear map $\mPhi_N: \ell_2^N \hookrightarrow \ell_1^{c\delta^{-2}N}$ such that 
\[
\left| \norm{\mPhi_N \vx}_1 - 1 \right| \leq \delta, \quad \forall \vx \in \mathbb{S}^{N-1}. 
\]
Similarly, $\ell_2^M$ is also embedded into $\ell_1^{c\delta^{-2}M}$ via $\mPhi_M: \ell_2^M \hookrightarrow \ell_1^{c\delta^{-2}M}$ so that \[
\left| \norm{\mPhi_M \vx}_1 - 1 \right| \leq \delta, \quad \forall \vx \in \mathbb{S}^{M-1}. 
\]

By the injectivity of the injective norm, $\mPsi = (\mI_K \otimes \mPhi_N) \otimes \mPhi_M$ embeds $\ell_1^K(\ell_2^N) \mathbin{\widecheck{\otimes}} \ell_2^M$ into $\ell_1^K(\ell_1^{c\delta^{-2}N}) \mathbin{\widecheck{\otimes}} \ell_1^{c\delta^{-2}M} \cong \ell_1^{c\delta^{-2}NK} \mathbin{\widecheck{\otimes}} \ell_1^{c\delta^{-2}M}$, where $\cong$ denotes the equivalence through an isometric isomorphism. 
Furthermore, it has been shown in Appendix~\ref{sec:ell1_to_ellinf} that $\ell_1^{c\delta^{-2}NK}$ (resp. $\ell_1^{c\delta^{-2}M}$) is isometrically embedded into $\ell_\infty^{2^{c\delta^{-2}NK}}$ (resp. $\ell_\infty^{2^{c\delta^{-2}M}}$). 
Therefore, due to the injectivity of the injective tensor norm, $\ell_1^{c\delta^{-2}NK} \mathbin{\widecheck{\otimes}} \ell_1^{c\delta^{-2}M}$ is embedded into $\ell_\infty^{2^{c\delta^{-2}NK}} \mathbin{\widecheck{\otimes}} \ell_\infty^{2^{c\delta^{-2}M}} \cong \ell_\infty^{2^{c\delta^{-2}(NK+M)}}$ via an isometric injection $\iota$. 
Moreover, the subspace $E = \mPsi(\ell_1^K(\ell_2^N) \mathbin{\widecheck{\otimes}} \ell_2^M)$ is also isometrically embedded to $\ell_\infty^{2^{c\delta^{-2}(NK+M)}}$. 

Let $Y = \ell_\infty^{2^{c\delta^{-2}(NK+M)}}$ and $F = \iota(E)$. 
Similar to Appendix~\ref{sec:ell1_to_ellinf}, by the Hahn-Banach theorem, there exists an isometric bijection $\varrho: Y^*/F^\perp \to F^*$ and the quotient map $q: Y^* \to Y^*/F^\perp$ is a metric surjection. 
Therefore, the map $\varrho \circ q: Y^* \to F^*$ is a metric surjection. 
Note that the restriction of $\iota^*$ on $F^*$, denoted by $\iota^*|_{F^*}: F^* \to E^*$, is an isometric bijection. Then the composition map $\mQ = \iota^*|_{F^*} \circ \varrho \circ q$ is a metric surjection. 
Moreover, the restriction of $\mPsi^*$ on $E^*$, denoted by $\mPsi^*|_{E^*}: E^* \to (\ell_\infty^K \mathbin{\widecheck{\otimes}} \ell_2^N) \mathbin{\widehat{\otimes}} \ell_2^M$ is a bijection. 
Therefore, there exists a map $(\mPsi^*|_{E^*})^{-1}$ such that $\mPsi^* \circ (\mPsi^*|_{E^*})^{-1}$ is the identity on $(\ell_\infty^K \mathbin{\widecheck{\otimes}} \ell_2^N) \mathbin{\widehat{\otimes}} \ell_2^M$. 
The embedding maps are illustrated in the following commutative diagram. 
\begin{center}
\begin{tikzcd}[row sep=large, column sep=large, ar symbol/.style = {draw=none,"#1" description,sloped}, isomorphic/.style = {ar symbol={\cong}}]
Y^* = \ell_1^{2^{c\delta^{-2}(NK+M)}} \arrow[d, two heads, "\mQ"] \arrow{rd}{\mT\mPsi^*\mQ} \\
E^* \arrow{r}{\mT \mPsi^*} & \ell_\infty^m(\ell_2^d) \\
(\ell_\infty^K \mathbin{\widecheck{\otimes}} \ell_2^N) \mathbin{\widehat{\otimes}} \ell_2^M \arrow[u, "(\mPsi^*|_{E^*})^{-1}"] \arrow[ru, "\mT"] 
\end{tikzcd}
\end{center}
Then the assertion follows from Lemma~\ref{mm} due to the surjectivity of the entropy number \cite[p. 12]{carl2008entropy} and the fact that $\norm{\mPsi} \leq (1+\delta)^2$.

\subsection{Proof of Lemma~\ref{lem:sup_quad}}
\label{sec:proof:lem:sup_quad}
Let $\vxi \in \R^{LMN}$ be a random vector defined by
\begin{equation}
\label{eq:xi}
    \vxi := \begin{bmatrix} \mathrm{vec}(\mB_{1,1})\\
    \mathrm{vec}(\mB_{2,1})\\
    \vdots \\
    \mathrm{vec}(\mB_{L,K})
    \end{bmatrix}.
\end{equation}
Since $\mathrm{vec}(\mB_{l,k})$'s are i.i.d. following $\mathcal{N}(\vzero, L^{-1} \mI_{MN})$, it follows that $\vxi \sim \mathcal{N}(\vzero, \mI_{LMN})$. 
Next, we define a matrix $\mQ_{\mZ} \in \R^{LK \times LMN}$ determined by $\mZ$ by
\begin{align*}
& \mQ_{\mZ} := \\
& \hspace{-10pt} \begin{bmatrix}
\mI_L \otimes \mathrm{vec}(\mZ_1)^* & 0& \cdots& 0 \\
0 & \mI_L \otimes\mathrm{vec}(\mZ_2)^* &  & 0\\
\vdots &  &\ddots & \\
0 & 0 & &\mI_L \otimes\mathrm{vec}(\mZ_K)^*
\end{bmatrix},
\end{align*}
where $\mZ_k = \mZ (\ve_k \otimes \mI_N)$ for $k \in [K]$ so that 
$\mZ = [\mZ_1 ~ \mZ_2 ~ \cdots ~ \mZ_K]$.
Then we have
\[
\sum_{l=1}^L \sum_{k=1}^K \langle \mA_{l,k}, \mZ \rangle^2 
= \norm{\mQ_{\mZ} \xi}_2^2.
\]
Furthermore, by taking the expectation on both sides, we obtain
\[
\E \sum_{l=1}^L \sum_{k=1}^K \langle \mA_{l,k}, \mZ \rangle^2 
= \E \norm{\mQ_{\mZ} \xi}_2^2
= \norm{\mQ_{\mZ}}_\mathrm{F}^2 
= L \norm{\mZ}_\mathrm{F}^2.
\]
Therefore, the left-hand side of \eqref{eq:lem:sup_quad} is written as the supremum of a second-order chaos as follows:
\begin{align*}
& \sup_{\mZ \in \kappa(\alpha,\beta)}
\left| L \norm{\mZ}_\mathrm{F}^2 - \sum_{l=1}^L \sum_{k=1}^K \langle \mA_{l,k}, \mZ \rangle^2 \right| \\
& = \sup_{\mZ \in \kappa(\alpha,\beta)}
\left| \norm{\mQ_{\mZ} \vxi}_2^2 - \E \norm{\mQ_{\mZ} \vxi}_2^2 \right|.
\end{align*}
We compute a tail bound on the right-hand side by using the results on suprema of chaos processes \cite{krahmer2014suprema}, which is summarized as Theorem~\ref{thm:kmr} in Appendix~\ref{sec:concentration}. 
To invoke Theorem~\ref{thm:kmr} for $\Delta = \{\mQ_{\mZ} : \mZ \in \kappa(\alpha,\beta)\}$, we derive upper bounds on the radii and the $\gamma_2$-functional of $\Delta$. 
The radii of $\Delta$ with respect to the spectral norm and the Frobenious norm satisfy
\[
d_\mathrm{S}(\Delta)
= \sup_{\mZ \in \kappa(\alpha, \beta)} \norm{\mQ_{\mZ}} 
= \sup_{\mZ \in \kappa(\alpha, \beta)} \max_{k \in [K]} \norm{\mZ_k}_\mathrm{F} 
\leq \alpha
\]
and
\[
d_\mathrm{F}(\Delta)
= \sup_{\mZ \in \kappa(\alpha, \beta)} \norm{\mQ_{\mZ}} 
= \sup_{\mZ \in \kappa(\alpha, \beta)} \sqrt{L} \norm{\mZ}_\mathrm{F} 
\leq \sqrt{LK} \alpha. 
\]
Note that $\norm{\mQ_{\mZ} - \mQ_{\mZ'}} = \norm{\mZ - \mZ'}_{\infty,\mathrm{F}}$ and $\kappa(\alpha,\beta) = \alpha B_{\infty,\mathrm{F}} \cap \beta B_\$$.
Then Dudley's inequality implies that the $\gamma_2$-functional of $\Delta$ is upper-bounded by
\begin{align*}
\gamma_2(\Delta) 
& \lesssim \int_0^\infty \sqrt{\ln N(\Delta, \eta B_\mathrm{S}) d\eta} \\
& \leq \beta \int_0^\infty \sqrt{\ln N\left(B_\$, \eta B_{\infty,\mathrm{F}}\right) d\eta}.
\end{align*}
Furthermore, Lemma~\ref{lemma:pi2a_nu1} implies 
\[
B_\$ 
\subset 
\sqrt{2N}
B_{(\ell_\infty^K \mathbin{\widecheck{\otimes}} \ell_2^N) \mathbin{\widehat{\otimes}} \ell_2^M}.
\]
Therefore, we obtain
\begin{align*}
& \int_0^\infty \sqrt{\ln N\left(B_\$, \eta B_{\infty,\mathrm{F}}\right)} d\eta \\
& \leq 
\int_0^\infty \sqrt{\ln N\left(\sqrt{N} B_{(\ell_\infty^K \mathbin{\widecheck{\otimes}} \ell_2^N) \mathbin{\widehat{\otimes}} \ell_2^M}, \eta B_{\infty,\mathrm{F}}\right)} d\eta \\
& \leq \sqrt{N} \int_0^\infty \sqrt{\ln N\left(B_{(\ell_\infty^K \mathbin{\widecheck{\otimes}} \ell_2^N) \mathbin{\widehat{\otimes}} \ell_2^M}, \eta B_{\ell_\infty^K(\ell_2^{MN})}\right)} d\eta \\
& \lesssim \sqrt{N} \sqrt{NK+M} (\ln K)^{3/2},
\end{align*}
where the last step follows from Lemma~\ref{lem:entropy_integral} together with the fact that 
\[
\norm{\mathrm{id}: (\ell_\infty^K \mathbin{\widecheck{\otimes}} \ell_2^N) \mathbin{\widehat{\otimes}} \ell_2^M \to \ell_\infty^K \mathbin{\widecheck{\otimes}} \ell_2^{MN}}_\mathrm{op} \leq 1, 
\]
which holds by Lemma~\ref{lem:interlacing}.
Combining these results provides
\[
\gamma_2(\Delta) \lesssim 
\beta \sqrt{N (NK+M) (\ln K)^3}.
\]
Then, the parameters $E$, $V$, and $U$ in Theorem~\ref{thm:kmr} are upper-bounded by
\begin{align*}
E {} & = \gamma_2(\Delta)
\left[ \gamma_2(\Delta) + d_{\mathrm{F}}(\Delta) \right] \\
& \lesssim \alpha^2 \sqrt{
(\beta/\alpha)^2 N (NK+M) (\ln K)^3} \\
& \quad \cdot \left( \sqrt{
(\beta/\alpha)^2 N (NK+M) (\ln K)^3} + \sqrt{LK} \right), \\
& = \alpha^2 L K \rho \left( \rho + 1 \right), \\
V {} & = d_{\mathrm{S}}(\Delta) \left[ \gamma_2(\Delta) + d_{\mathrm{F}}(\Delta) \right] \\
& \lesssim \alpha^2 \left( \sqrt{
(\beta/\alpha)^2 N (NK+M) (\ln K)^3} + \sqrt{LK} \right) \\
& = \alpha^2 \sqrt{LK} (\rho + 1), \\
U {} & = d_{\mathrm{S}}^2(\Delta) \leq \alpha^2.
\end{align*}
Then, by plugging in these parameters into Theorem~\ref{thm:kmr}, we obtain that 
\begin{equation}
\label{eq:lem:sup_quad_tmp}
\begin{aligned}
& \sup_{\mZ \in \kappa(\alpha,\beta)}
\left| \norm{\mZ}_\mathrm{F}^2 - \sum_{l=1}^L \sum_{k=1}^K \langle \mA_{l,k}, \mZ \rangle^2 \right| \\
& \lesssim 
\alpha^2 \ln(2\zeta^{-1}) + \alpha^2 L K \left( \rho + 1 \right) \left(\rho + \sqrt{\frac{\ln(2\zeta^{-1})}{LK}} \right),
\end{aligned}
\end{equation}
holds with probability $1-\zeta$. 
Finally, by choosing $C$ in \eqref{eq:sample_rate} large enough, we have $\rho \leq 1$, which further simplifies \eqref{eq:lem:sup_quad_tmp} into \eqref{eq:lem:sup_quad}. 
This completes the proof. 

\subsection{Proof of Lemma~\ref{lemma:ub_cor}}
\label{sec:proof:lemma:ub_cor}
Let $\tnorm{\cdot}$ denote the norm defined so that the unit norm ball is $\kappa(\alpha,\beta)$, i.e. 
\[
\kappa(\alpha,\beta) = \{ \mX \in \mathbb{R}^{M \times NK} : \tnorm{\mX} \leq 1 \}.
\]
Then the left-hand side of \eqref{eq:lemma_ub_cor} is written as
\begin{equation}
\label{eq:lemma_ub_cor2}
\sup_{\mZ \in \kappa(\alpha, \beta)} \sum_{l=1}^L \sum_{k=1}^K \langle \mA_{l,k}, \mZ \rangle w_{l,k} 
= \tnorm{\sum_{l=1}^L \sum_{k=1}^K \langle \mA_{l,k}, \mZ \rangle w_{l,k}}_*,  
\end{equation}
where $\tnorm{\cdot}_*$ denotes the dual norm of $\tnorm{\cdot}$. 
Conditioned on $\mA_{l,k}$'s, the quantity on the right-hand side of \eqref{eq:lemma_ub_cor2} becomes a Gaussian empirical process. 
Due to \cite[Theorem~4.7]{pisier1999volume}, it holds with probability $1-\frac{\zeta}{3}$ that 
\begin{equation}
\label{eq:ub_noise_dualnorm}
\begin{aligned}
& \tnorm{\sum_{l=1}^L \sum_{k=1}^K \langle \mA_{l,k}, \mZ \rangle w_{l,k}}_* \\
& \leq \sigma \E_{(g_{l,k})} \tnorm{\sum_{l=1}^L \sum_{k=1}^K \langle \mA_{l,k}, \mZ \rangle g_{l,k}}_* \\
& + 
\sigma \pi 
\sqrt{ \frac{\ln(6\zeta^{-1})}{2} \cdot \sup_{\vert\vert\vert \mZ \vert\vert\vert \leq 1} \sum_{l=1}^L \sum_{k=1}^K \langle \mA_{l,k}, \mZ \rangle^2},
\end{aligned}
\end{equation}
where $g_{l,k}$'s are i.i.d. Gaussian with zero mean and unit variance. 

The last term in the right-hand side of \eqref{eq:ub_noise_dualnorm} is upper-bounded by using the following result. 
Due to Lemma~\ref{lem:sup_quad}, there exists a numerical constant $C$, for which it holds with probability $1-\frac{\zeta}{3}$ that
\begin{align*}
& \sum_{l=1}^L \sum_{k=1}^K \langle \mA_{l,k}, \mZ \rangle^2 
\leq \norm{\mZ}_\mathrm{F}^2 \\
& \quad + C \left[ \frac{\alpha^2 \ln(6\zeta^{-1})}{L} + \alpha^2 K \left(\rho + \sqrt{\frac{\ln(6\zeta^{-1})}{LK}} \right) \right]
\end{align*}
for all $\mZ \in \kappa(\alpha,\beta)$. 
Furthermore, we also have
\[
\sup_{\mZ \in \kappa(\alpha,\beta)} \norm{\mZ}_\mathrm{F}^2 
\leq \sup_{\mZ \in \kappa(\alpha,\beta)} K \norm{\mZ}_{\infty,\mathrm{F}}
\leq \alpha^2 K.
\]
Therefore, we obtain
\begin{align*}
& \sup_{\mZ \in \kappa(\alpha,\beta)} \sum_{l=1}^L \sum_{k=1}^K \langle \mA_{l,k}, \mZ \rangle^2 \\
& \lesssim 
\frac{\alpha^2 \ln(6\zeta^{-1})}{L} 
+ \alpha^2 K \left(\rho + \sqrt{\frac{\ln(6\zeta^{-1})}{LK}} + 1 \right).    
\end{align*}

Furthermore, due to \citep[Equation (4.9)]{ledoux2013probability}, the expectation term in the right-hand side of \eqref{eq:ub_noise_dualnorm} is upper-bounded by
\begin{align*}
& \E_{(g_{l,k})} \tnorm{\sum_{l=1}^L \sum_{k=1}^K \langle \mA_{l,k}, \mZ \rangle g_{l,k}}_* \\
& \lesssim \sqrt{\ln(LK+1)}
\E_{(\epsilon_{l,k})} \tnorm{\sum_{l=1}^L \sum_{k=1}^K \langle \mA_{l,k}, \mZ \rangle \epsilon_{l,k}}_*,    
\end{align*}
where $(\epsilon_{l,k})$ is a Rademacher sequence, i.e. $\epsilon_{l,k}$'s are independent copies of random variable $\epsilon$ satisfying $\P{\epsilon = 1} = \P{\epsilon = -1} = \frac{1}{2}$.

Then, due to the symmetry of the distribution of $\mA_{l,k}$'s, we obtain the following identity, which holds in the sense of distribution with respect to $\mA_{l,k}$'s: 
\begin{samepage}
\begin{align}
& \E_{(\epsilon_{l,k})} \tnorm{\sum_{l=1}^L \sum_{k=1}^K \langle \mA_{l,k}, \mZ \rangle \epsilon_{l,k}}_* \\
& = \sup_{\vert\vert\vert \mZ \vert\vert\vert \leq 1} \left| \sum_{l=1}^L \sum_{k=1}^K \langle \epsilon_{l,k} \mA_{l,k}, \mZ \rangle \right| \nonumber \\ 
& = \sup_{\vert\vert\vert \mZ \vert\vert\vert \leq 1} \left| \sum_{l=1}^L \sum_{k=1}^K \langle \mA_{l,k}, \mZ \rangle \right| \nonumber \\
& = \sup_{\vert\vert\vert \mZ \vert\vert\vert \leq 1} \sum_{l=1}^L \sum_{k=1}^K \langle \mA_{l,k}, \mZ \rangle \nonumber \\
& = \sup_{\vert\vert\vert \mZ \vert\vert\vert \leq 1} \sum_{l=1}^L \sum_{k=1}^K \langle \ve_k^* \otimes \mB_{l,k}, \mZ \rangle, \label{eq:rademacher_complexity}
\end{align}
where the third step follows due to the symmetry in $\kappa(\alpha,\beta)$ and the last step used $\mA_{l,k} = \ve_k^* \otimes \mB_{l,k}$. \end{samepage}

Let $\vxi \in \R^{LMN}$ be defined in \eqref{eq:xi}. 
Furthermore, with a shorthand notation $\mZ_k := \mZ (\ve_k \otimes \mI_N) \in \R^{M \times N}$, we define a column vector $\vf_{\mZ} \in \R^{LMN}$ given by
\begin{equation*}
\vf_{\mZ} := 
\begin{bmatrix}
\vct{1}_{L,1} \otimes \mathrm{vec}(\mZ_1) \\
\vct{1}_{L,1} \otimes \mathrm{vec}(\mZ_2) \\
\vdots \\
\vct{1}_{L,1} \otimes \mathrm{vec}(\mZ_N) \\
\end{bmatrix},
\end{equation*}
where $\vct{1}_{L,1} \in \R^L$ denotes the column vector with all entries set to $1$. 
Then the last term in \eqref{eq:rademacher_complexity} is written as
\[
\sup_{\vert\vert\vert \mZ \vert\vert\vert \leq 1} \sum_{l=1}^L \sum_{k=1}^K \langle \mB_{l,k}, \mZ (\ve_k \otimes \mI_N) \rangle 
= \sup_{\vert\vert\vert \mZ \vert\vert\vert \leq 1} \langle \vf_{\mZ}, \bm{\xi} \rangle.
\]
Note that the right-hand side is the supremum of a Gaussian process. 
To obtain an upper bound, we will use Lemma~\ref{lemma:gaussianmax} in Appendix~\ref{sec:concentration}. 
To invoke Lemma~\ref{lemma:gaussianmax} for the set $\Delta = \{ \vf_{\mZ} : \tnorm{\mZ} \leq 1 \}$ and $\bm{\xi} \sim \mathcal{N}(\vzero,\mI_{LMN})$, we compute the diameter and covering number of $\Delta$ with respect to the $\ell_2$-norm. 
Since 
\begin{align*}
\norm{\vf_{\mZ} - \vf_{\mZ'}}_2 
= \sqrt{L} \norm{\mZ - \mZ'}_\mathrm{F}
\leq \sqrt{LK} \norm{\mZ - \mZ'}_{\infty,\mathrm{F}}, 
\end{align*}
it follows that $\mathrm{diam}(\Delta) \leq \sqrt{LK} \alpha$ and
\[
N(\Delta, \eta B_2) 
\leq 
N\left( B_{\vert\vert\vert\cdot\vert\vert\vert}, \eta (LK)^{-1/2} B_{\infty,\mathrm{F}} \right).
\]
Furthermore, since
\begin{align*}
B_{\vert\vert\vert\cdot\vert\vert\vert} 
& = \kappa(\alpha,\beta) 
= \alpha B_{\infty,\mathrm{F}} \cap \beta B_\$,
\end{align*}
by Lemma~\ref{lemma:pi2a_nu1}, we have 
\begin{align*}
B_{\vert\vert\vert\cdot\vert\vert\vert} 
& \subset \beta B_\$ 
\subset \beta \sqrt{2N} B_{(\ell_\infty^K \mathbin{\widecheck{\otimes}} \ell_2^N) \mathbin{\widehat{\otimes}} \ell_2^M}.    
\end{align*}
Therefore, we obtain
\begin{align*}
& \int_0^\infty \sqrt{\ln N\left( B_{\vert\vert\vert\cdot\vert\vert\vert}, \eta (LK)^{-1/2} B_{\infty,\mathrm{F}} \right)} d\eta \\
& \leq 
\beta \sqrt{LKN} \int_0^\infty \sqrt{\ln N\left(B_{(\ell_\infty^K \mathbin{\widecheck{\otimes}} \ell_2^N) \mathbin{\widehat{\otimes}} \ell_2^M}, \eta B_{\ell_\infty^K(\ell_2^{MN})}\right)} d\eta \\
& \lesssim
\beta \sqrt{LKN} \sqrt{NK+M} (\ln K)^{3/2},
\end{align*}
where the last step follows from Lemma~\ref{lem:entropy_integral} and Lemma~\ref{lem:interlacing}.
By plugging in this result to Lemma~\ref{lemma:gaussianmax}, we obtain that 
\begin{align*}
& \sup_{\vert\vert\vert \mZ \vert\vert\vert \leq 1} \langle \vf_{\mZ}, \bm{\xi} \rangle \\
& \lesssim 
\beta \sqrt{LKN} \sqrt{NK+M} (\ln K)^{3/2} + \alpha \sqrt{LK \ln(3\zeta^{-1})} \\
& = \alpha LK \rho + \alpha \sqrt{LK \ln(3\zeta^{-1})}
\end{align*}
holds with probability $1-\frac{\zeta}{3}$. 

\subsection{Proof of Theorem~\ref{thm:minimax}}

We establish the minimax lower bound in Theorem~\ref{thm:minimax} by following the two-step strategy outlined below. 
We first show that there exists a packing set of $\kappa(\alpha, \beta)$ of a desirable size and a packing density. 
Then a minimax bound is derived via a multi-way hypothesis testing argument and Fano's inequality. 

Let us first recall the notion of a packing set (e.g. see \citep[Definition~4.2.4]{vershynin2018high}). 
A subset $\mathcal{P}$ of a metric space $\mathcal{S}$ is called $\epsilon$-packing of $\mathcal{S}$ if $d(x, x') > \epsilon$ for all distinct $x,x' \in \mathcal{P}$, where the parameter $\epsilon$ denotes the packing density. 
The following lemma constructs a packing set of $\kappa(\alpha, \beta)$ with respect to the metric induced by the Frobenius norm. 

\begin{lemma} 
Let $\gamma \leq 1$ satisfy that $\frac{\beta^2}{\gamma^2 \alpha^2}$ is an integer. 
Then there exists a subset $\mathcal{H} \subset \kappa(\alpha, \beta)$ with cardinality 
\[ 
|\mathcal{H}| = \left\lfloor \exp \left( \frac{(\beta/\alpha)^2(NK \vee M)}{16 \gamma^2} \right) \right\rfloor 
\] 
with the following properties:
\begin{enumerate}
    \item Every $\mH \in \mathcal{H}$ satisfies that $\mathrm{rank}(\mH) \leq \frac{\beta^2}{\alpha^2 \gamma^2}$ and each entry is from $\{ \pm \frac{\gamma \alpha}{\sqrt{MN}} \}$, thereby $\norm{\mH}_{\infty, F} = \gamma \alpha$ and $\norm{\mH}_\mathrm{F}^2 = K \gamma^2 \alpha^2$.

    \item Any two distinct $\mH^i, \mH^j \in \mathcal{H}$ satisfy $$\norm{\mH^i - \mH^j}_\mathrm{F}^2 
    \geq \frac{K\gamma^2 \alpha^2}{2}.$$
\end{enumerate}
\label{lm:pack_set}
\end{lemma}

\begin{IEEEproof}
We adapt the proof of \cite[Lemma 3.1]{cai2016matrix} to our setting. 
The idea is to show the existence of a packing set by the empirical method. 
We first consider the case where $NK \geq M$. 
Let $S = \lfloor \exp(\frac{(\beta/\alpha)^2 KN}{16 \gamma^2}) \rfloor$ and $B = \frac{\beta^2}{\alpha^2\gamma^2}$. 
We generate $\mH^1,\dots,\mH^S$ as independent copies of a random matrix $\mH$ constricted as follows. 
The entries of the first $B$ rows of $\mH$ are i.i.d. following the uniform distribution on $\{ \pm \frac{\gamma \alpha}{\sqrt{MN}} \}$. 
The remaining rows are determined from the first $B$ rows by
\begin{align*}
H_{m,n} = H_{m',n}, \quad & \forall m, m' \in [M] : m' \equiv m \, (\mathrm{mod} \, B), ~ \\
& \forall n \in [NK]. 
\end{align*}
Since the magnitude of all entries of $\mH$ are fixed to the constant $\frac{\gamma\alpha}{\sqrt{MN}}$, it follows that $\norm{\mH}_{\infty,\mathrm{F}} = \gamma \alpha$ and $\norm{\mH}_\mathrm{F} = \sqrt{K} \gamma \alpha$. Furthermore, by Lemma~\ref{lem:interlacing}, we also have
\[
\norm{\mH}_\$ \leq \sqrt{B} \norm{\mH}_{\infty,\mathrm{F}} = \frac{\beta}{\gamma \alpha} \cdot \gamma \alpha = \beta,
\]
thereby, $\mH^i \in \kappa(\alpha,\beta)$ for all $i \in [S]$, or equivalently, $\mathcal{H} \subset \kappa(\alpha,\beta)$. 

For any $\mH^i \neq \mH^j$, we have 
\begin{align*}
    \norm{\mH^i - \mH^j}_\mathrm{F}^2 
    & = \sum_{m=1}^M \sum_{n=1}^{NK} (H^i_{m,n} - H^j_{m,n})^2 \\
    & \geq \left\lfloor \frac{M}{B} \right\rfloor \sum_{m=1}^{B} \sum_{n=1}^{NK} (H^i_{m,n} - H^j_{m,n})^2 \\ 
    & \geq \frac{4\alpha^2 \gamma^2}{MN} \left\lfloor \frac{M}{B} \right\rfloor \sum_{m=1}^{B} \sum_{n=1} ^{KN} \delta_{m,n},
\end{align*} 
where $\delta_{m,n}$'s are i.i.d. symmetric Bernoulli random variables. 
By Hoeffding's inequality, we obtain 
\begin{equation*}
    \mathbb{P} \left ( \sum_{m=1}^B \sum_{n=1}^{KN} \delta_{m,n} \leq \frac{BNK}{4} \right ) 
    \leq e^{-\frac{BKN}{8}}.
\end{equation*} 
By the union bound argument over all $\binom{S}{2}$ possible distinct pairs $(\mH^i,\mH^j)$, we obtain that 
\begin{equation*}
    \min_{i \neq j} \norm{\mH^i - \mH^j}_\mathrm{F}^2 
    > \alpha^2 \gamma^2 \left\lfloor \frac{M}{B} \right\rfloor \frac{BKN}{MN} 
    \geq \frac{\alpha^2 \gamma^2  K}{2}
\end{equation*} 
holds with probability at least $1 - \binom{S}{2} \exp(-\frac{BNK}{8}) \geq 1 - \frac{S^2}{2} \exp(-\frac{BNK}{8}) \geq \frac{1}{2}$. 
In other words, the second property is satisfied with nonzero probability, thereby, there exists such an instance. 
If $M > NK$, then we construct $\mH^*$ by the same procedure. Then the existence of a desired packing set is shown similarly. 
This concludes the proof. 
\end{IEEEproof}

\begin{lemma}[Equivalence to multiple hypothesis testing][Lemma 6.2, \cite{lee2021approximately}] 
\label{lem:mht}
Let $\calH$ be a $\delta$-packing set of $\kappa (\alpha, \beta)$ and let $\widetilde{\mH} = \mathop{\mathrm{argmin}}_{\mH \in \calH}\|\mH - \widehat{\mH} \|_F$. 
Then we have 
\begin{equation*}
    \inf_{\widehat{\mH}} \sup_{\mH \in \kappa(\alpha,\beta)} \E \|\widehat{\mH} - \mH\|_\mathrm{F}^2 
    \geq 
    \frac{\delta^2}{4} \min_{\widetilde{\mH} \in \mathcal{H}}\P { \widetilde{\mH} \neq \mH^* }, 
\end{equation*} 
where $\mH^*$ is uniformly distributed over $\mathcal{H}$. 

\end{lemma}

\noindent We now proceed to a lower bound on $\min_{\widetilde{\mH} \in \calH} \P{ \widetilde{\mH} \neq \mH^* }$. 
To this end, we use the following version of Fano's inequality stated in \cite{cai2016matrix}. 
\begin{lemma}[Fano's inequality] 
\label{lem:fano}
Let $\widetilde{\mH} = \mathop{\mathrm{argmin}}_{\mH \in \calH}\|\mH - \widehat{\mH} \|_F$.
Then we have  
\begin{equation}
\label{eq:Fanos}
\begin{aligned}
& \mathbb{P} (\widetilde{\mH} \neq \mH^*) \\
& \geq 1 - \frac{ {\binom{|\mathcal{H}|}{2}} ^{-1} \sum_{i \neq j} \mathbb{E}_{(\mB_{l,i})} D_\mathrm{KL} (\mH^i \,\|\, \mH^j)+ \ln 2}{\ln | \mathcal{H}|}\,,
\end{aligned}
\end{equation}
where $D_\mathrm{KL}(\mH^i \,\|\, \mH^j)$ denotes the Kullback–Leibler divergence between the joint distributions of $y_{l,k}$'s in the measurement model \eqref{eq:sensing_model} conditioned on measurement matrices $\mB_{l,k}$'s for $\mH^i$ and $\mH^j$. 
\end{lemma}

\noindent It remains to compute the KL divergence in \eqref{eq:Fanos} so that we can invoke Fano's inequality in Lemma~\ref{lem:fano}.
The joint probability density of $y_{l,k}$'s given $\mB_{l,k}$'s is given by 
\begin{align*}
& p\left( \{ y_{l,k} \} | \{\mB_{l,k} \} \right) \\
& = \prod_{l = 1}^{L} \prod_{k = 1}^{K} \frac{1}{\sqrt{2 \pi \sigma^2}} \exp\left(- \frac{\left( y_{l,k} - \langle \mB_{l,k}, \mH_i^{k} \rangle \right )^2}{2L\sigma^2}\right).
\end{align*}
Then we obtain
\begin{align*}
& \ln \left ( \frac{ p(\vy | \mH_i) }{ p(\vy | \mH_j) } \right ) \\
&= \sum_{l=1}^L \sum_{k=1}^K \frac{(y_{l,i} - \langle \mB_{l,k}, \mH_k^j \rangle)^2 - (y_{l,k} - \langle \mB_{l,k}, \mH_k^i \rangle)^2}{ 2\sigma^2} \\
&= \sum_{l=1}^L \sum_{k=1}^K \frac{\langle \mB_{l,k}, \mH_k^i-\mH_k^j\rangle^2}{ 2\sigma^2} \\
&\quad + \sum_{l=1}^L \sum_{k=1}^K \frac{(y_{l,k} - \langle \mB_{l,k}, \mH_k^i \rangle) \langle \mB_{l,k}, \mH_k^i-\mH_k^j\rangle}{ 2\sigma^2},
\end{align*}
where $\mH_k^i = \mH^i (\ve_k \otimes \mI_N)$ denotes the $k$th block of $\mH^i$ of size $M \times N$ for $k \in [K]$. 
Hence it follows that 
\begin{align*}
& D_\mathrm{KL}(\mH_i \,\|\, \mH_j) 
= \int_{-\infty}^\infty p( \vy | \mH_i ) \ln \left ( \frac{ p(\vy | \mH_i) }{ p(\vy | \mH_j) } \right ) d y \\
& = \frac{1}{2\sigma^2} \sum_{l=1}^L \sum_{k=1}^K \langle \mB_{l,k} , \mH_i^k - \mH_j^k\rangle^2.
\end{align*} 
Furthermore, by taking the expectation with respect to $\mB_{l,k}$'s, we obtain
\begin{equation}
\mathbb{E}_{(\mB_{l,k})} D_\text{KL} (\mH^i \,\|\, \mH^j) = \frac{L}{2\sigma^2}\norm{\mH^i - \mH^j}_\mathrm{F}^2
\label{eq:KLdiv}
\end{equation}

We consolidate the above sequence of results to establish the minimax lower bound in Theorem~\ref{thm:minimax}. 
Recall that the packing set $\mathcal{H}$ by Lemma~\ref{lm:pack_set} satisfies $\norm{\mH^i}_\mathrm{F} = \sqrt{K} \gamma \alpha$ for all $\mH^i \in \mathcal{H}$. It immediately follows that $\norm{\mH^i - \mH^j}_\mathrm{F}^2 \leq 4 K \gamma^2 \alpha^2$ for all $\mH^i, \mH^j \in \mathcal{H}$.
Furtheremore, we have
\[
\ln |\mathcal{H}| 
\leq \frac{(\beta/\alpha)^2(NK \vee M)}{16 \gamma^2}. 
\]
Plugging in these result together with \eqref{eq:KLdiv} into \eqref{eq:Fanos}, we obtain
\begin{align*}
& \mathbb{P} (\widetilde{\mH} \neq \mH^*) \\
& \geq 1 - \frac{16 \gamma^2}{(\beta/\alpha)^2 (NK \vee M)} \left( \frac{2 L K \gamma^2 \alpha^2}{\sigma^2} + \ln 2 \right)
\geq \frac{1}{2}
\end{align*} 
provided that $\gamma^4 \leq \frac{(\beta/\alpha)^2(NK \vee M) \sigma^2}{128 L K \alpha^2}$ and $(\beta/\alpha)^2 (NK \vee M) \geq 48$.

If it is satisfied that $\frac{(\beta/\alpha)^2(NK \vee M) \sigma^2}{128 L K \alpha^2} \geq 1$, then we choose $\gamma = 1$. 
In this case, by Lemma~\ref{lem:mht}, we obtain
\begin{align*}
     \inf_{\widehat{\mH}} \sup_{\mH \in \kappa(\alpha,\beta)} \E \|\widehat{\mH} - \mH\|_\mathrm{F}^2 
     & \geq \frac{\delta^2}{4} \cdot \frac{1}{2} 
     \geq \frac{K \alpha^2}{16},
\end{align*} 
where the second inequality follows since the packing density of $\mathcal{H}$ was $\delta = \gamma \alpha \sqrt{\frac{K}{2}}$. 
This implies 
\begin{align*}
     \inf_{\widehat{\mH}} \sup_{\mH \in \kappa(\alpha,R)} \E \frac{1}{K} \|\widehat{\mH} - \mH\|_\mathrm{F}^2 & 
     \geq \frac{\alpha^2}{16}.
\end{align*} 

\noindent Otherwise, we choose $\gamma = \left( \frac{(\beta/\alpha)^2(NK \vee M) \sigma^2}{128 L K \alpha^2} \right)^{1/4}$ so that 
\begin{align*}
& \inf_{\widehat{\mH}} \sup_{\mH \in \kappa(\alpha,R)} \E \|\widehat{\mH} - \mH\|_\mathrm{F}^2 \\
& \geq \frac{\delta^2}{4} \cdot \frac{1}{2} \geq \frac{K \alpha \sigma}{16} \cdot \sqrt{\frac{(\beta/\alpha)^2(NK \vee M)}{128LK}},
\end{align*} 
which follows from Lemma~\ref{lem:mht}.
Therefore, we have
\begin{align*}
& \inf_{\widehat{\mH}} \sup_{\mH \in \kappa(\alpha,\beta)} \E \frac{1}{K} \|\widehat{\mH} - \mH\|_\mathrm{F}^2 \\
& \geq \frac{\alpha^2}{16\cdot 8\sqrt{2}} \cdot \frac{\sigma}{\alpha}  \sqrt{\frac{(\beta/\alpha)^2(NK \vee M)}{LK}}.
\end{align*} 
Finally, combining the two results, we obtain
\begin{align*}
& {\inf_{\widehat{\mH}}} {\sup_{\mH\in \kappa(\alpha, R)}} \frac{1}{K} \E \|\widehat{\mH} - \mH\|_{F}^2 \\
& \geq \frac{\alpha^2}{16} \left(1 \wedge \frac{\sigma}{8\sqrt{2} \alpha} \sqrt{\frac{(\beta/\alpha)^2(NK \vee M)}{LK}} \right ).
\end{align*}
This completes the proof.

\subsection{ADMM algorithms}
\label{sec:admm}
The optimization formulations in \eqref{eq:compute_dollarnorm} and \eqref{eq:estimator} can be rewritten into a standard semidefinite program and be solved by off-the-shelf solvers like SeDuMi \cite{sturm1999using} or SDPT3 \cite{toh1999sdpt3}. 
However, these software packages do not scale well to large instances. To alleviate the limitation, we develop \textit{Alternating Direction Method of Multipliers} (ADMM) algorithms, wherein each subproblem admits a closed-form solution or casts as a simple program easily solved by standard linear algebra packages.

\subsubsection{ADMM algorithm to compute the $\$$-norm in \eqref{eq:compute_dollarnorm}}

We first rewrite the optimization formulation in \eqref{eq:compute_dollarnorm} into an equivalent problem with a set of auxiliary variables:
\begin{equation}
\label{eq:equivalent_compute_dollarnorm}
\def\arraystretch{1.25}
\begin{array}{lcl}
& \displaystyle \mathop{\mathrm{minimize}}_{\beta,\mW_1,\mW_2,\mE} ~ & \displaystyle \beta \\
& \mathrm{subject~to} &\displaystyle \mathrm{trace}(\mW_{1}) \leq \beta \\
& & \displaystyle 
\mathrm{trace}\left((\ve_k^* \otimes \mI_N) \mW_2 (\ve_k \otimes \mI_N)\right) \leq \beta, \\ 
& & \hspace{4.35cm} \forall k \in [K] \\
& &  \mE = \displaystyle \begin{bmatrix}
\mW_{1} & \mX \\
\mX^* & \mW_{2}
\end{bmatrix} \\
& & \mE \succeq \mathbf{0}. 
\end{array}
\end{equation}
Then an augmented Lagrangian function of \eqref{eq:equivalent_compute_dollarnorm} is obtained by penalizing the equality constraints as
\begin{equation*}
\begin{aligned}
& \mathcal{L}_{\rho}\left(\beta,\mW_1,\mW_2,\mE,\mPhi\right) \\
& = \beta+\langle\mPhi,\mE- \begin{bmatrix}
\mW_1 & \mX \\
\mX^* & \mW_2
\end{bmatrix}\rangle 
+\frac{\rho}{2} \, \left\|\mE-\begin{bmatrix}
\mW_1 & \mX \\
\mX^* & \mW_2
\end{bmatrix}\right\|_{\mathrm{F}}^2,
\end{aligned}
\end{equation*} 
where $\mPhi\in\R^{(M+KN)\times(M+KN)}$ denotes a dual variable. 
ADMM finds a global minimizer to the convex program in \eqref{eq:equivalent_compute_dollarnorm} by minimizing $\mathcal{L}_{\rho}$ with respect to each of the primal variables $\beta,\mW_1,\mW_2,\mE$ sequentially followed by the gradient ascent update of the dual variable $\mPhi$ \cite[Section~3.1]{boyd2011distributed}. 

For brevity, we introduce the following shorthand notations. 
We decompose $\mE, \mPhi \in \R^{(M+NK) \times (M+NK)}$ into four blocks as
\begin{equation}
\label{eq:simple_notation}
\begin{aligned}
&\mE:= \begin{bmatrix}
\mE_{11} & \mE_{12} \\
\mE_{12}^* & \mE_{22}
\end{bmatrix} 
\quad \text{and} \quad 
\mPhi =\begin{bmatrix}
\mPhi_{11} & \mPhi_{12} \\
\mPhi_{12}^* & \mPhi_{22}
\end{bmatrix},
\end{aligned}
\end{equation}
where the size of each block is given by $\mE_{11}, \mPhi_{11} \in \R^{M \times M}$, $\mE_{12}, \mPhi_{12} \in \R^{M \times NK}$, and $\mE_{22}, \mPhi_{22} \in \R^{NK \times NK}$. 
Further, the $k$th block of size $M \times N$ of $\mX \in \R^{M \times NK}$ given by $\mX (\ve_k \otimes \mI_N)$ is denoted by $\mX_k$. 
Similarly, the $k$th diagonal block of size $N \times N$ of $\mW_2 \in \R^{NK \times NK}$ given by $(\ve_k^* \otimes \mI_N) \mW_2 (\ve_k \otimes \mI_N)$ is denoted by $\mW_{2,k}$. 

Given the above shorthand notations, we describe the update rules of the ADMM algorithm. 
First, we consider the updates of the first block of primal variables. 
We update $\beta$, $\mW_1$, and $\{\mW_{2,k}\}_{k=1}^K$ by solving the following optimization problem: 
\begin{equation}
\label{eq:subproblem1}
\def\arraystretch{1.25}
\begin{array}{ll} 
\displaystyle \mathop{\mathrm{minimize}}_{\mW_1,\{\mW_{2,k}\}_{k=1}^K,\beta} & \displaystyle \frac{\rho}{2} \norm{\mW_1 - \mA}_F^2 + \frac{\rho}{2} \sum_{k=1}^K \norm{\mW_{2,k} - \mB_k}_F^2 + \beta \\
\text{subject~to} & \mathrm{trace}(\mW_1) \leq \beta \\
& \mathrm{trace}(\mW_{2,k}) \leq \beta, \quad \forall k \in [K], 
\end{array}
\end{equation} where 
\[
\mA:=\mE_{11}+{\rho}^{-1}\mPhi_{11}, \quad \text{and} \quad \mB:={\mE}_{22}+{\rho}^{-1}\mPhi_{22}.
\]
For fixed $\beta$, minimization decouples over the other variables and the optimal solution is given by 
\begin{align*}
\widehat{\mW}_1
= \mathop{\mathrm{argmin}}_{\trace({\mW_{1}})\leq\beta} \left\|\mW_1-\mA\right\|_{\mathrm{F}}^2
\end{align*}
and
\begin{align*}
\widehat{\mW}_{2,k}
= \mathop{\mathrm{argmin}}_{\trace({\mW_{2}})\leq\beta} \left\|\mW_{2,k}-\mB_k\right\|_{\mathrm{F}}^2.
\end{align*} 
Furthermore, $\widehat{\mW}_1$ and $\widehat{\mW}_{2,k}$ are expressed in a closed-form respectively given by
\begin{equation}
\label{eq:update_W1}
\widehat{\mW}_1 = \mA - \left(\frac{\max(\mathrm{trace}(\mA)-
\beta,0)}{M}\right) \mI_M
\end{equation} 
and
\begin{equation}
\label{eq:update_W2k}
\widehat{\mW}_{2,k} = \mB_k - \left(\frac{\max(\mathrm{trace}(\mB_k)-
\beta,0)}{M}\right) \mI_N.
\end{equation} 
By plugging in the expression of the optimal solutions in \eqref{eq:update_W1} and \eqref{eq:update_W2k} for fixed $\beta$ into \eqref{eq:subproblem1}, the optimization formulation in \eqref{eq:subproblem1} reduces to the minimization of a univariate function given by
\begin{equation}
\label{eq:minimizebeta}
\begin{aligned}
f(\beta) & = \beta+\frac{\rho \max(\mathrm{trace}(\mA)-\beta,0)^2}{2M} \\ 
& +\sum_{k=1}^K\frac{\rho \max(\mathrm{trace}(\mB_k)-\beta,0)^2}{2N}. 
\end{aligned}
\end{equation}
Due to the monotonicity of the summands in the right-hand side of \eqref{eq:minimizebeta}, the global minimizer $\hat{\beta}$ can be found by the bisection search on the interval from $0$ to $\max(\mathrm{trace}(\mA),\max_{k\in[K]}\mathrm{trace}(\mB_k))$. 
Once $\beta$ is updated as $\hat{\beta}$, then $\mW_1$ (resp. $\mW_{2,k}$) will be updated as $\widehat{\mW}_1$ by \eqref{eq:update_W1} (resp. $\widehat{\mW}_{2,k}$ by \eqref{eq:update_W2k}). 
The off-diagonal blocks of $\mW_2$ are updated by
\begin{align*}
& (\ve_j^* \otimes \mI_N) \mW_2 (\ve_k \otimes \mI_N) \\
& = (\ve_j^* \otimes \mI_N) (\mE_{22}+\rho^{-1}\mPhi_{22}) (\ve_k \otimes \mI_N), 
\quad j \neq k \in [K]. 
\end{align*}
Next, the second block of primal variables consists of $\mE$, which is updated as the solution to 
\begin{align*}
\widehat{\mE}
& = \mathop{\mathrm{argmin}}_{\mE\succeq \vzero} \left\langle\mPhi,\mE\right\rangle+\frac{\rho}{2} \, \left\|\mE-\begin{bmatrix}
\mW_{1} & \mX \\
\mX^* & \mW_{2}
\end{bmatrix}\right\|_{\mathrm{F}}^2 \\
& = \mathcal{P}_{\mathbb{S}_+^{M+KN}} \left(\begin{bmatrix}
\mW_{1} & \mX \\
\mX^* & \mW_{2}
\end{bmatrix}-\rho^{-1}\mPhi\right),
\end{align*}
where $\mathbb{S}_+^{M+KN}$ denotes the cone of positive semidefinite matrices of size $(M+NK)$. 
Finally, the dual variable $\mPhi$ is updated by gradient ascent with step size $\rho$. 

\subsubsection{ADMM algorithm for the convex estimator in \eqref{eq:estimator}}
The optimization in \eqref{eq:estimator_sdp} is equivalently reformulated with an auxiliary variable as 
\begin{equation}
\label{eq:eqivalent_problem}
\def\arraystretch{1.25}
\begin{array}{lcl}
& \displaystyle \mathop{\mathrm{minimize}}_{\mX, \mZ, \mW_1, \mW_2} ~ & \displaystyle \sum_{l=1}^L \sum_{k=1}^K \left( y_{l,k} - \langle \mB_{l,k}, \mX (\ve_k \otimes \mI_N) \rangle  \right)^2 \\
& \mathrm{subject~to} &\displaystyle \mathrm{trace}(\mW_{1}) \leq \beta \\
& & \displaystyle 
\mathrm{trace}\left((\ve_k^* \otimes \mI_N) \mW_2 (\ve_k \otimes \mI_N)\right) \leq \beta, \\ 
& & \hspace{4.55cm} k \in [K], \\
& & \left\|\mX\right\|_{\infty,\mathrm{F}}\leq\alpha\\
& & \mZ = \displaystyle \begin{bmatrix}
\mW_{1} & \mX \\
\mX^* & \mW_{2}
\end{bmatrix} \\
& & \mZ \succeq \mathbf{0}. 
\end{array}
\end{equation}
An augmented Lagrangian function is written as
\begin{equation*}
\begin{aligned}
& \mathcal{L}_{\rho}\left(\mX,\mW_1,\mW_2,\mZ,\mPsi\right) \\
& =\sum_{l=1}^L \sum_{k=1}^K \left( y_{l,k} - \langle \mB_{l,k}, \mX (\ve_k \otimes \mI_N) \rangle  \right)^2 \\ 
& +\langle\mPsi,\mZ- \begin{bmatrix}
\mW_1 & \mX \\
\mX^* & \mW_2
\end{bmatrix}\rangle 
+\frac{\rho}{2}\left\|\mZ-\begin{bmatrix}
\mW_1 & \mX \\
\mX^* & \mW_2
\end{bmatrix}\right\|_{\mathrm{F}}^2,
\end{aligned}
\end{equation*} 
where $\mPsi\in\R^{(M+KN)\times(M+KN)}$ denotes a dual variable. 
Then the ADMM algorithm iterates the minimization of $\mathcal{L}_{\rho}$ with respect to primal-variable blocks $(\mX,\mW_1,\mW_2)$ and $\mZ$ followed by the gradient ascent update of the dual variable $\mPsi$ as shown below. 
For brevity, we consider the decomposition of $\mZ, \mPsi \in \R^{(M+NK) \times (M+NK)}$ into four blocks given by
\begin{align*}
&\mZ:= \begin{bmatrix}
\mZ_{11} & \mZ_{12} \\
\mZ_{12}^* & \mZ_{22}
\end{bmatrix} 
\quad \text{and} \quad 
\mPsi =\begin{bmatrix}
\mPsi_{11} & \mPsi_{12} \\
\mPsi_{12}^* & \mPsi_{22}
\end{bmatrix},
\end{align*}
where the size of each block is given by $\mZ_{11}, \mPsi_{11} \in \R^{M \times M}$, $\mZ_{12}, \mPsi_{12}, \in \R^{M \times NK}$, and $\mZ_{22}, \mPsi_{22} \in \R^{NK \times NK}$. 

First, we consider the update of $\mX$, $\mW_1$, and $\mW_2$ in the first block. 
Note that the minimization of $\mathcal{L}_\rho$ only with respect to $\mX$ reduces to a norm-constrained least squares problem. Due to the blockwise structure in the measurement model, it decouples over blocks of $\mX$ as
\begin{equation}
\label{eq:update_Xk}
\begin{aligned}
\mathop{\mathrm{minimize}}_{\|\mX_k\|_{\mathrm{F}}\leq \alpha} &
\displaystyle \sum_{l=1}^L \left( y_{l,k} - \langle \mB_{l,k}, \mX_k \rangle  \right)^2 \\ 
& \quad - \left\langle\mPsi_{12,k},\mX_k\right\rangle+\frac{\rho}{2}\left\|\mZ_{12,k}-\mX_k\right\|_{\mathrm{F}}^2    
\end{aligned}
\end{equation}
for $k \in [K]$. Then, \eqref{eq:update_Xk} is written as 
\begin{equation}
\label{eq:update_Xk_normconstraint}
\begin{array}{cc}
\displaystyle \mathop{\mathrm{minimize}}_{\|\mathrm{vec}({\mX}_k)\|_2\leq \alpha} &
\displaystyle \left\|\vb-\mQ\mathrm{vec}({\mX}_k)\right\|_2^2
\end{array}
\end{equation} 
for $\mQ$ and $\vb$ satisfying 
$\mQ^*\mQ:=2\widetilde{\mB}_k^*\widetilde{\mB}_k +\rho\mI$ and $\mQ^*\vb:=2\widetilde{\mB}_k^*{\bf y}_k+\mathrm{vec}(\mPsi_{12,k})+\rho\, \mathrm{vec}(\mZ_{12,k})$ where
\[
    \widetilde{\mB}_k := \begin{bmatrix} \mathrm{vec}(\mB_{1,k})^*\\
    \mathrm{vec}(\mB_{2,k})^*\\
    \vdots \\
    \mathrm{vec}(\mB_{L,k})^*
    \end{bmatrix}
    \quad \text{and} \quad
    \vy_k:=\begin{bmatrix} y_{1,k}\\
    y_{2,k}\\
    \vdots \\
    y_{L,k}
    \end{bmatrix}.
\] 
Then \eqref{eq:update_Xk_normconstraint} becomes a norm-constrained least square problem.
Since \eqref{eq:update_Xk_normconstraint} satisfies the Slater's condition, the minimizer is obtained by the Karush–Kuhn–Tucker (KKT) conditions through the Lagrangian function 
\begin{equation}
\label{eq:update_Xk_Lagrangian}
\mathcal{L}(\mathrm{vec}({\mX}_k),\lambda):=\|\vb-\mQ\mathrm{vec}({\mX}_k)\|_2^2+\lambda\left(\left\|\mathrm{vec}({\mX}_k)\right\|_2^2-\alpha^2\right)
\end{equation}
given by 
\begin{equation}
\label{eq:KKT_conditions}
\begin{aligned}
& \|\mathrm{vec}({\mX}_k)\|_2\leq\alpha, \\
& \lambda\geq 0, \\ 
& \lambda\,(\|\mathrm{vec}({\mX}_k)\|_2-\alpha)=0, \\ 
& \left(\mQ^*\mQ-\mQ^*\vb\right)+\lambda\,\mathrm{vec}({\mX}_k)=0.    
\end{aligned}
\end{equation}
The optimal Lagrangian multiplier $\lambda^\star$ can be found by a binary search as outlined below. 
The unique minimizer to \eqref{eq:update_Xk_Lagrangian}, denoted by $\widehat{\mX}_k^{\lambda}$, is given by
\begin{equation}
\label{eq:Xksolution}
\begin{aligned}
\mathrm{vec}(\widehat{\mX}_k^{\lambda}) 
& =  \left(2\widetilde{\mB}_k^*\widetilde{\mB}_k +(\rho+\lambda)\mI\right)^{-1} \\ 
& \quad \cdot 
\left(2\widetilde{\mB}_k^*{\bf y}_k+\mathrm{vec}(\mPsi_{12,k})+\rho\, \mathrm{vec}(\mZ_{12,k})\right).
\end{aligned}
\end{equation}
Then the KKT conditions \eqref{eq:KKT_conditions} will be satisfied by the optimal Lagrange multiplier $\lambda^\star$ and $\mathrm{vec}(\widehat{\mX}_k^{\lambda^\star})$.  
Note that the solution in \eqref{eq:Xksolution} satisfies the last condition in \eqref{eq:KKT_conditions} for all $\lambda\geq0$. Furthermore, since $\|\mathrm{vec}(\widehat{\mX}_k^{\lambda})\|_2$ is a decreasing function of $\lambda$, the optimal $\lambda^\star$ can be found by a bisection method. 
Moreover, since $\widetilde{\mB}_k$ depends only on $\mB_{l,k}$'s, which do not vary over iterations, the solution in \eqref{eq:Xksolution} is easily obtained from a pre-compute the eigenvalue decomposition of $(\widetilde{\mB}_k^*\widetilde{\mB}_k)^{-1}$. 
The update of $\mW_1$ is given by
\begin{align*}
\widehat{\mW}_1
= \mathop{\mathrm{argmin}}_{\trace({\mW_{1}})\leq\beta} -\langle\mPsi_{11},\mW_{1}\rangle+\frac{\rho}{2}\left\|\mZ-\mW_{1}\right\|_{\mathrm{F}}^2,
\end{align*}
which yields a closed-form expression
\begin{align*}
\widehat{\mW}_1
& = \mZ_{11}+\rho^{-1}\mPsi_{11} \\ 
& - \left(\frac{\max(\mathrm{trace}(\mZ_{11}+\rho^{-1}\mPsi_{11})-
\beta,0)}{M}\right) \mI_M.
\end{align*}
Similarly, the diagonal blocks of $\mW_2$ are updated as
\begin{align*}
\widehat{\mW}_{2,k}
& = \mZ_{22,k}+\rho^{-1}\mPsi_{22,k} \\ 
& - \left(\frac{\max(\mathrm{trace}(\mZ_{22,k}+\rho^{-1}\mPsi_{22,k})-
\beta,0)}{N}\right) \mI_N.    
\end{align*}
The off-diagonal blocks of $\mW_2$ are copied from the corresponding blocks of $\mZ_{22}+\rho^{-1}\mPsi_{22}$.
Next, the primal variable $\mZ$ in the second block is updated by
\begin{align*}
\widehat{\mZ}
& = \mathop{\mathrm{argmin}}_{\mZ\succeq \vzero} \left\langle\mPsi,\mZ\right\rangle+\frac{\rho}{2} \, \left\|\mZ-\begin{bmatrix}
\mW_{1} & \mX \\
\mX^* & \mW_{2}
\end{bmatrix}\right\|_{\mathrm{F}}^2 \\
& = \mathcal{P}_{\mathbb{S}_+^{M+KN}} \left(\begin{bmatrix}
\mW_{1} & \mX \\
\mX^* & \mW_{2}
\end{bmatrix}-\rho^{-1}\mPsi\right),
\end{align*}
where $\mathbb{S}_+^{M+KN}$ denotes the cone of positive semidefinite matrices of size $(M+NK)$. 
Finally, the dual variable $\mPsi$ is updated by gradient ascent with step size $\rho$.  

For fast convergence, we adopt a varying step size for the dual ascent \cite[Section~3.4.1]{boyd2011distributed}, in which $\rho$ is updated in each iteration by keeping the primal and dual residual norms within a constant factor of each other. 
Furthermore, we employed a stopping criterion based on the feasibility and relative change of primal variables \cite[Section~3.3.1]{boyd2011distributed}, 
which has been widely used in practice.

\end{document}